\documentclass{amsart}[12pt]

\usepackage{epigraph}
\usepackage{amssymb}
\usepackage{stmaryrd}
\usepackage[all]{xy}
\usepackage{mathrsfs}
\usepackage{color}
\usepackage{geometry}
\usepackage{enumitem}
\usepackage{color}

\numberwithin{itemcounter}{subsection}

\def\k{\mathbf{k}}
\def\C{\mathbb{C}}
\def\d{\mathbf{d}}
\def\H{\mathbf{H}}
\def\N{\mathbb{N}}
\def\Z{\mathbb{Z}}
\def\Q{\mathbb{Q}}
\def\dim{\text{dim}}
\def\Hom{\text{Hom}}
\def\Ext{\text{Ext}}
\def\g{\mathfrak{g}}

\def\e{\mathbf{e}}
\def\n{\mathbf{n}}

\def\v{\mathbf{v}}
\def\w{\mathbf{w}}

\theoremstyle{plain}

\newtheorem{theorem}{Theorem}[section]

\newtheorem{lemma-definition}[theorem]{Lemma-Definition}
\newtheorem{definition-lemma}[theorem]{Definition-Lemma}

\newtheorem{proposition}[theorem]{Proposition}
\newtheorem{conjecture}[theorem]{Conjecture}

\theoremstyle{definition}

\theoremstyle{remark}

\numberwithin{equation}{section}

\author[O.~Schiffmann]{Olivier Schiffmann}
\address[Olivier Schiffmann]{D\'epartement de Math\'ematiques, Universit\'e de Paris-Sud Paris-Saclay, B\^at. 307, 91405 Orsay Cedex, France}
\curraddr{}
\email{olivier.schiffmann@math.u-psud.fr}

\title[Kac polynomials and Lie algebras associated to quivers and curves]{Kac polynomials and Lie algebras associated to \\ quivers and curves}

\begin{document}

\maketitle

\epigraph{\emph{Ah! si loin des carquois, des torches et des fl\`eches,\newline On se sauvait un peu vers des choses... plus fra\^iches !}}
{ E. Rostand, Cyrano de Bergerac, Acte III, Sc\`ene 7 }

\vspace{.1in}

\section{Kac polynomials for quivers and curves}

\vspace{.1in}

\paragraph{\textbf{1.1. Quivers.}} Let $Q$ be a locally finite quiver with vertex set $I$ and edge set $\Omega$. For any finite field $\mathbb{F}_q$ and any dimension vector $\d \in \N^I$, let $A_{Q,\d}(\mathbb{F}_q)$ be the number of absolutely indecomposable representations of $Q$ over $\mathbb{F}_q$, of dimension $\d$. Kac proved the following beautiful result~:

\begin{theorem}[Kac, \cite{Kac81}] There exists a (unique) polynomial $A_{Q,\d}(t) \in \Z[t]$ such that for any finite field $\mathbb{F}_q$, 
$$A_{Q,\d}(\mathbb{F}_q)=A_{Q,\d}(q).$$
Moreover $A_{Q,\d}(t)$ is independent of the orientation of $Q$ and is monic of degree $1-\langle \d, \d \rangle$.
\end{theorem}

 Here $\langle\;,\;\rangle$ is the Euler form, see Section~\textbf{2.1.} The fact that the number of (absolutely) indecomposable $\mathbb{F}_q$-representations behaves polynomially in $q$ is very remarkable : the  absolutely indecomposable representations only form a \textit{constructible} substack of the stack $\mathcal{M}_{Q,\d}$ of representations of $Q$ of dimension $\d$ and we count them here \textit{up to isomorphism}, i.e. \textit{without} the usual orbifold measure. Let us briefly sketch the idea of a proof. Standard Galois cohomology arguments ensure that it is enough to prove that the number of indecomposable $\mathbb{F}_q$-representations is given by a polynomial $I_{Q,\d}$ in $q$ (this polynomial is not as well behaved as or as interesting as $A_{Q,\d}$). Next, by the Krull-Schmidt theorem, it is enough to prove that the number of \textit{all} $\mathbb{F}_q$-representations of dimension $\d$ is itself given by a polynomial in $\mathbb{F}_q$. This amounts to computing the (orbifold) volume of the \textit{inertia} stack $\mathcal{IM}_{Q,\d}$ of $\mathcal{M}_{Q,\d}$; performing a unipotent reduction in this context, we are left to computing the volume of the stack
$\mathcal{N}il_{Q,\d}$ parametrizing pairs $(M,\phi)$ with $M$ a representation of dimension $\d$ and $\phi \in \text{End}(M)$ being nilpotent. Finally, we use a Jordan stratification of $\mathcal{N}il_{Q,\d}$ and easily compute the volume of each strata in terms of the volumes of the stacks $\mathcal{M}_{Q,\d'}$ for all $\d'$. This actually yields an explicit formula for $A_{Q,\d}(t)$ (or $I_{Q,\d}(t)$), see \cite{Hua}.  We stress that beyond the case of a few quivers (i.e. those of finite or affine Dynkin type) it is unimaginable to classify and construct all indecomposable representations; nevertheless, the above theorem says that we can \textit{count} them.

\vspace{.05in}

The positivity of $A_{Q,\d}(t)$ was only recently\footnote{the special case of an indivisible dimension vector was proved earlier by Crawley-Boevey and Van den Bergh, see \cite{CBVdB}} proved~:

\begin{theorem}[Hausel-Letellier-Rodriguez-Villegas, \cite{HLRVAnnals}]\label{T:positivity} For any $Q$ and any $\d$ we have $A_{Q,\d}(t) \in \N[t]$.
\end{theorem}

The desire to understand the meaning of the Kac polynomials $A_{Q,\d}(t)$, for instance as the dimension of certain natural graded vector spaces (such as the cohomology of some algebraic variety), has been a tremendous source of inspiration in geometric representation theory; it has lead to the development over the years of a beautiful and very rich
theory relating various moduli spaces of representations of quivers to Lie algebras and quantum groups, yielding geometric constructions of a host of important objects in representation theory (such as canonical and crystal bases, Yangians, quantum affine algebras, their highest weight or finite-dimensional representations, etc.). We will present here some recent developments in this area, in particular in relation to the structure and representation theory of Kac-Moody algebras and \textit{graded Borcherds algebras}.

\vspace{.1in}

\paragraph{\textbf{1.2.~Curves.}} Now let $X$ be a smooth projective curve defined over some finite field $\mathbb{F}_q$. For any $r \geq 0, d \in \mathbb{Z}$ let $A_{r,d}(X)$ be the number of isomorphism classes of absolutely indecomposable coherent sheaves $\mathcal{F}$ over $X$ of rank $r$ and degree $d$. It turns out that this number may also be computed explicitly and depends 'polynomially' on $X$. In order to make sense of this, we need to introduce a few notations. The action of the geometric Frobenius $Fr$ on the $l$-adic cohomology group $H_{et}^1(X_{\overline{\mathbb{F}_q}},\overline{\mathbb{Q}_l})$ is semisimple with eigenvalues $\{\sigma_1, \ldots, \sigma_{2g}\}$ which may be ordered so that $\sigma_{2i-1}\sigma_{2i}=q$ for all $i$. Moreover, $Fr$ belongs to the general symplectic group $GSp(H_{et}^1(X_{\overline{\mathbb{F}_q}},\overline{\mathbb{Q}_l}))$ relative to the intersection form on $H_{et}^1(X_{\overline{\mathbb{F}_q}},\overline{\mathbb{Q}_l})$. We may canonically identify the character ring
of $GSp(H_{et}^1(X_{\overline{\mathbb{F}_q}},\overline{\mathbb{Q}_l}))$ with $R_g= \overline{\mathbb{Q}_l}[T_g]^{W_g}$ where
$$T_g=\{(\eta_1, \ldots, \eta_{2g}) \in \mathbb{G}_m^{2g}\;|\; \eta_{2i-1}\eta_{2i}=\eta_{2j-1}\eta_{2j}\;\forall\;i,j\}, \qquad W_g = (\mathfrak{S}_2)^g \rtimes \mathfrak{S}_g$$
are the maximal torus, resp. Weyl group of $GSp(2g,\overline{\mathbb{Q}_l})$. We can evaluate any element $f \in R_g$ on $Fr_X$ for any smooth projective curve
$X$ of genus $g$, defined over some\footnote{technically, of characteristic different from $l$} finite field $\mathbb{F}_q$; concretely, we have
$$f(Fr_X)=f(\sigma_1, \sigma_2, \ldots, \sigma_{2g}).$$
We will say that a quantity depending on $X$ is \textit{polynomial} if it is the evaluation of some element $f \in R_g$. For instance the size of the field $\mathbb{F}_q$ of definition of $X$ is polynomial in $X$ (take $f=\eta_{2i-1}\eta_{2i}$ for any $i$), as is the number $|X(\mathbb{F}_q)|$ of $\mathbb{F}_q$-rational points of $X$ (take $f=1-\sum_i \eta_i + q$); this is also true of symmetric powers $S^lX$ of $X$. We extend the above definition to the case $g=0$ by setting $R_g=\overline{\mathbb{Q}_l}[q^{\pm 1}]$.

\vspace{.1in}

\begin{theorem}[S., \cite{Sannals}]\label{T:Kaccurve} For any fixed $g,r,d$ there exists a (unique) polynomial $A_{g,r,d} \in R_g$ such that for any smooth projective curve $X$ of genus $g$ defined over a finite field, the number of absolutely indecomposable coherent sheaves on $X$ of rank $r$ and degree $d$ is equal to $A_{g,r,d}(Fr_X)$. Moreover, $A_{g,r}$ is monic of leading term $q^{1+(g-1)r^2}$.
\end{theorem}

The finiteness of the number of absolutely indecomposable coherent sheaves of fixed rank and degree is a consequence of Harder-Narasimhan reduction theory; it suffices to observe that any sufficiently unstable coherent sheaf is decomposable as its Harder-Narasimhan filtration splits in some place. Let us say a few words about the proof of Theorem~\ref{T:Kaccurve}~: like for quivers, standard Galois cohomology arguments combined with the Krull-Schmidt theorem reduce the problem to counting \textit{all} isomorphism classes of coherent sheaves of rank $r$ and degree $d$; unfortunately, this number is infinite as soon as $r >0$ hence it is necessary to introduce a suitable truncation of the category
$\text{Coh}(X)$; there are several possibilities here, one of them being to consider the category $\text{Coh}^{\geq 0}(X)$ of positive coherent sheaves, i.e. sheaves whose HN factors all have positive degree. One is thus lead to compute the volumes of the inertia stacks $\mathcal{I}\mathcal{M}_{X,r,d}^{\geq 0}$ of the stacks $\mathcal{M}_{X,r,d}^{\geq 0}$ positive sheaves of rank $r$ and degree $d$ on $X$. Finally, after a unipotent and Jordan reduction similar to the case of quivers, this last computation boils down to the evaluation of the integral of certain Eisenstein series over the truncated stacks $\mathcal{M}_{X,r,d}^{\geq 0}$; this is performed using a variant of Harder's method for computing the volume of $\mathcal{M}_{X,r,d}$, see \cite{Sannals} for details. The proof is constructive and yields an explicit but complicated formula for $A_{g,r,d}$. Hence, just like for quivers, although classifying all indecomposable vector bundles for curves of genus $g >1$ is a wild problem, it is nevertheless possible to count them. The explicit formula for $A_{g,r,d}$ was later combinatorially very much simplified by Mellit, see \cite{Mellit17a}, who in particular proved the following result~:

\begin{theorem}[Mellit, \cite{Mellit17a}] The polynomial $A_{g,r,d}$ is independent of $d$.
\end{theorem}

Thanks to this theorem, we may simply write $A_{g,r}$ for $A_{g,r,d}$. Like their quivery cousins, Kac polynomials of curves satisfy some positivity and integrality property, namely

\begin{theorem}[S., \cite{Sannals}] For any $g, r$, $A_{g,r} \in \text{Im}(\mathbb{N}[-\eta_1, \ldots, -\eta_g]^{W_g} \to R_g).$
\end{theorem}

In other words, $A_{g,r}$ is a $W_g$-symmetric polynomial in $-\eta_1, \ldots ,-\eta_{2g}$ with positive integral coefficients. However, this is \textit{not} the most natural form of positivity~:

\begin{conjecture}\label{C:kacpolcat} For any $g,r$ there exists a canonical (\emph{non}-virtual) finite-dimensional representation $\mathbb{A}_{g,r}$ of $GSp(2g,\overline{\mathbb{Q}_l})$ such that
$A_{g,r}=\tau (ch(\mathbb{A}_{g,r}))$, where $\tau \in Aut(R_g)$ is the involution mapping $\eta_i$ to $-\eta_i$ for all $i$.
\end{conjecture}

The paper \cite{Hua} contains many examples of Kac polynomials for quivers. Let us give here the examples of Kac polynomials for curves (of any genus) for $r=0,1,2$; we will write $q$ for $\eta_{2i-1}\eta_{2i}$.
$$A_{g,0}=1-\sum_i \eta_i + q,\; \qquad A_{g,1}=\prod_{i=1}^{2g} (1-\eta_i),$$
$$A_{g,2}=\prod_{i=1}^{2g}(1-\eta_i) \cdot \left( \frac{\prod_{i} (1-q\eta_i)}{(q-1)(q^2-1)} -\frac{\prod_i (1+\eta_i)}{4(1+q)} + \frac{\prod_i (1-\eta_i)}{2(q-1)} \left[ \frac{1}{2}-\frac{1}{q-1}-\sum_i \frac{1}{1-\eta_i}\right]\right).$$
For $r=0,1$ we recognize the number of $\mathbb{F}_q$-rational points of $X$ and $Pic^0(X)$ respectively. In particular, we have $\mathbb{A}_{g,0}=\mathbb{C} \oplus \mathbb{V} \oplus det(V)^{1/g}$ and $\mathbb{A}_{g,1}=\Lambda^\bullet V$, where $V$ is the standard $2g$-dimensional representation of $GSp(2g,\overline{\mathbb{Q}_l})$.

\vspace{.1in}

Motivated by the analogy with quivers, it is natural to try to seek a representation-theoretic meaning to the Kac polynomials $A_{g,r}$. What is the analog, in this context, of the Kac-Moody algebra associated to a quiver~? Are the Kac polynomials related to the Poincar\'e polynomials of some interesting moduli spaces ?
Although it is still much less developed than in the context of quivers, we will illustrate this second theme through two examples of applications of the theory of \textit{Hall algebras} of curves~: a case of geometric Langlands duality (in the neighborhood of the trivial local system) and the computation of the Poincar\'e polynomial of the moduli spaces of semistable Higgs bundles on smooth projective curves.

\vspace{.1in}

\paragraph{\textbf{1.3.~Quivers vs.~Curves}} Kac polynomials of quivers and curves are not merely related by an analogy : they are connected through the following observation, which comes by comparing the explicit formulas for $A_{Q,\d}$ and $A_{g,r}$. Let $S_g$ be the quiver with one vertex and $g$ loops.

\begin{proposition}[Rodriguez-Villegas]\label{P:strangerelation} For any $r$ we have $A_{S_g,r}(1)=A_{g,r}(0,\ldots, 0)$.\end{proposition}

This relation between the constant term of $A_{g,r}$ and the sum of all the coefficients of $A_{S_g,r}$ has a very beautiful conjectural conceptual explanation
in terms of the mixed Hodge structure of the (twisted) genus $g$ character variety, we refer the interested reader to \cite{HRVInvent}. We will provide another 
conceptual explanation in terms of the geometric Langlands duality in Section~\textbf{5.4}.

\vspace{.1in}

\noindent
\textit{Remark.} There is an entirely similar story for the category of coherent sheaves on a smooth projective curve equipped with a (quasi-)parabolic structure along an effective divisor $D$, see \cite{SDukeI}, \cite{Lin}, \cite{Mellit17b}. Proposition~\ref{P:strangerelation} still holds in this case, with the quiver $S_g$ being replaced by a quiver with a central vertex carrying $g$ loops, to which are attached finitely many type $A$ branches, one for each point in $D$.

\vspace{.1in}

\noindent
\textit{Plan of the paper.} In sections~\textbf{2, 3} and \textbf{4} we describe various Lie-theoretical (\textbf{2.2, 2.5, 3.3, 4.4, 4.6.}) or geometric (\textbf{4.1, 4.4, 4.5, 4.6.}) incarnations (some conjectural) of the Kac polynomials for quivers; in particular we advocate the study of a certain \textit{graded} Borcherds algebra $\widetilde{\g}_Q$ canonically associated to a quiver, which is a graded extension of the usual Kac-Moody algebra $\g_Q$ attached to $Q$. From Section~\textbf{5} onward, we turn our attention to curves and present several algebraic or geometric constructions suggesting the existence of some hidden, combinatorial and
Lie-theoretical structures controlling such things as dimensions of spaces of cuspidal functions (\textbf{6.2.}), or Poincar\'e polynomials of moduli spaces of (stable) Higgs bundles on curves (\textbf{7.2, 7.3}).

\medskip

\section{Quivers, Kac polynomials and Kac-Moody algebras}

\medskip

We begin with some recollections of some classical results in the theory of quivers and Kac-Moody algebras, including the definitions of Hall algebras and Lusztig nilpotent varieties, and their relations to (quantized) enveloping algebra of Kac-Moody algebras. This is related to the \textit{constant term} $A_{Q,\d}(0)$ of the Kac polynomials $A_{Q,\d}(t)$. The theory is classically developed for quivers without edge loops. The general case, important for applications, is more recent. 

\vspace{.1in}

\paragraph{\textbf{2.1. Kac-Moody algebras from quivers.}} Let $Q=(I,\Omega)$ be a finite quiver without edge loops. For any field $\k$ we denote by $\text{Rep}_\k(Q)$ the abelian category of $\k$-representations of $Q$ and for any dimension vector $\d \in \N^I$ we denote by $\mathcal{M}_{Q,\d}$ the stack of $\d$-dimensional representations of $Q$. The category $\text{Rep}_\k(Q)$ is of global dimension at most one (and exactly one if $\Omega$ is non-empty) with finite-dimensional Ext spaces. As a result, the stack $\mathcal{M}_Q=\bigsqcup_{\d} \mathcal{M}_{Q,\d}$ is smooth. The first relation to Lie theory arises when considering the Euler forms
\begin{align}
\langle M,N\rangle = \dim\;\Hom(M,N) - \dim\; \Ext^1(M,N), \qquad (M,N)=\langle M,N \rangle + \langle N,M\rangle.
\end{align}
Let $c_{ij}$ be the number of arrows in $Q$ from $i$ to $j$, and let $C=(c_{ij})_{i,j \in I}$ be the adjacency matrix. Set $A=2Id - C - {}^tC$.
The Euler forms $\langle\;,\;\rangle$ and $(\,,\,)$ factor through the map
$$\underline{\dim}~: K_0(\text{Rep}_{\k}Q) \to \mathbb{Z}^I, \qquad M \mapsto (\dim\; M_i)_i$$
and are given by
$$\langle \d, \d' \rangle= {}^t\d \,(Id -C) \,\d', \qquad (\d,\d')={}^t\d \,A \,\d'.$$
Note that $\dim\;\mathcal{M}_{Q,\d}=-\langle \d,\d\rangle$. 
Now, $A$ is a (symmetric) generalized Cartan matrix in the sense of \cite{KacBook}, to which is attached a Kac-Moody algebra $\g_Q$. The Euler lattice $(\Z^I, (\,,\,))$ of $\text{Rep}_{\k}Q$ is identified to the root lattice\footnote{we apologize for the unfortunate --yet unavoidable-- clash between the traditional notations for the root lattice and for the quiver} $Q_{\g_Q}$ of $\g_Q$ together with its standard Cartan pairing, via the map $\Z^I \to Q_{\g_Q}, \d \mapsto \sum_i d_i \alpha_i$ (here $\alpha_i$ are the simple roots of $\g_Q$). Accordingly, we denote by $\g_Q=\bigoplus_{\d}\g_Q[\d]$ the root space decomposition of $\g_Q$.

\vspace{.1in} 

\paragraph{\textbf{2.2. Kac's theorem and the constant term conjecture.}} Kac proved that $A_{Q,\d}(t) \neq 0$ if and only if $\d$ belongs to the root system $\Delta^+$ of $\mathfrak{g}_Q$; this generalizes the famous theorem of Gabriel \cite{Gabriel} which concerns the case of finite-dimensional $\mathfrak{g}_Q$. Moreover, he made the following conjecture, which was later proved by Hausel (see \cite{CBVdB} for the case of indivisible $\d$)~:

\begin{theorem}[Hausel, \cite{HauselKac}]\label{T:constant} For any $\d \in \Delta^+$, we have $A_{Q,\d}(0) = \dim\; \g_{Q}[\d]$.
\end{theorem}

We will sketch a proof of Theorems~\ref{T:positivity} and \ref{T:constant} using cohomological Hall algebras in Section~\textbf{4.4}. These proofs are different from, but related to the original proofs of Hausel et al.

\vspace{.1in}

In the remainder of Section~\textbf{2}, we discuss
two (related) constructions of the Kac-Moody algebra $\mathfrak{g}_Q$~: as the \textit{spherical Hall algebra} of the categories $\text{Rep}_\k(Q)$ for $\k$ a finite field, or in terms of the complex geometry of the \textit{Lusztig nilpotent varieties} $\Lambda_{Q,\d}$.

\vspace{.1in}

\paragraph{\textbf{2.3. Ringel-Hall algebras}} It is natural at this point to ask for an actual construction of $\g_Q$ using the moduli stacks $\mathcal{M}_{Q,\d}$ of representations of $Q$.
This was achieved by Ringel and Green, \cite{Ringel}, \cite{Green} using the \textit{Hall} (or \textit{Ringel-Hall}) algebra of $\text{Rep}_{\k}Q$. Let $\k$ be again a finite field and set  $\mathcal{M}_Q=\bigsqcup_{\d} \mathcal{M}_{Q,\d}$. The set of $\k$-points $\mathcal{M}_{Q}(\k)=\bigsqcup_{\d} \mathcal{M}_{Q,\d}(\k)$ is by construction the set of $\k$-representations of $Q$ (up to isomorphism). Put
$$\H_Q=\bigoplus_{\d} \H_Q[\d], \qquad \H_Q[\d]:=\text{Fun}(\mathcal{M}_{Q,\d}(\k) ,\C)$$
and consider the following convolution diagram
\begin{equation}\label{E:conv}
\xymatrix{ \mathcal{M}_Q \times \mathcal{M}_Q & \widetilde{\mathcal{M}}_Q \ar[l]_-q \ar[r]^-p & \mathcal{M}_Q}
\end{equation}
where $\widetilde{\mathcal{M}}_Q$ is the stack parametrizing short exact sequences 
$$\xymatrix{0 \ar[r] & M \ar[r] & R \ar[r] & N \ar[r]& 0}$$
in $\text{Rep}_{\k}Q$ (alternatively, inclusions $M \subset R$ in $\text{Rep}_{\k}Q$); the map $p$ assigns to a short exact sequence as above its middle term $R$; the map $q$ assigns to it its end terms $(N,M)$. The stack $\widetilde{\mathcal{M}}_Q$ can be seen as parametrizing extensions between objects in $\text{Rep}_\k Q$. The map $p$ is proper and, because $\text{Rep}_{\k}Q$ is of homological dimension one, the map $q$ is a stack vector bundle whose restriction to $\mathcal{M}_{Q,\d} \times \mathcal{M}_{Q,\d'}$ is of rank
$-\langle \d',\d\rangle$. Following Ringel and Green, we put $v=(\#\k)^{\frac{1}{2}}$, let $\mathbf{K}=\C[k_i^{\pm 1}]_{i}$ be the group algebra of $\Z^I$ and we equip $\widetilde{\H}_Q:=\H_Q \otimes \mathbf{K}$ with the structure of a $\N^I$-graded bialgebra by setting
\begin{equation}\label{E:prodhall}
f \cdot g =v^{\langle \d, \d'\rangle} p_!q^* (f \otimes g), \qquad k_\n f k_\n^{-1}=v^{(\d,\n)} f 
\end{equation}
\begin{equation}\label{E:prodhall2}
\Delta (f)=\sum_{\d'+\d''=\d} v^{\langle \d',\d''\rangle} (q_{\d',\d''})_!p_{\d',\d''}^*(f)  \qquad \Delta (k_i)=k_i \otimes k_i
\end{equation}
for $f \in \H_Q[\d], \;g \in \H_{Q}[\d'], \n \in \Z^I$ and $i \in I$. The bialgebra $\widetilde{\H}_Q$ is equipped with the nondegenerate Hopf pairing 
\begin{equation}\label{E:prodhall3}(f{k}_\d,g{k}_{\d'})=v^{(\d,\d')}\int_{\mathcal{M}_Q} f \overline{g}, \qquad \forall\; f,g \in \H_Q
\end{equation}
 (see \cite{SLectures} for more details and references).  Denoting by 
$\{\epsilon_i\}_i$ the canonical basis fo $\Z^I$ we have $\mathcal{M}_{\epsilon_i} \simeq B\mathbb{G}_m$ for all $i$ since
$\mathcal{M}_{\epsilon_i}$ has as unique object the simple object $S_i$ of dimension $\epsilon_i$. 
We define the \textit{spherical} Hall algebra $\widetilde{\H}^{\text{sph}}_Q$ as the subalgebra of $\widetilde{\H}_Q$ generated by $\mathbf{K}$ and the elements $1_{S_i}$ for $i \in I$. It is a sub-bialgebra. 

\begin{theorem}[Ringel, \cite{Ringel}, Green, \cite{Green}] The assignement $E_i \mapsto 1_{S_i}, K^{\pm 1}_i \mapsto k^{\pm 1}_i$ for $i \in I$ induces an isomorphism $\Psi~:U_v(\mathfrak{b}_Q) \stackrel{\sim}{\to} \widetilde{\H}^{\text{sph}}_Q$ between the positive Borel subalgebra of the Drinfeld-Jimbo quantum enveloping algebra $U_v(\g_Q)$ and the spherical Hall algebra of $Q$.
\end{theorem}
{\text{sph}}
One may recover the entire quantum group $U_v(\g_Q)$ as the (reduced) Drinfeld double $\mathbf{D}\widetilde{\H}_Q^{\text{sph}}$ of $\widetilde{\H}_Q$ (see e.g. \cite[Lec.\;5]{SLectures}).

\vspace{.1in}

\paragraph{\textbf{2.4.\;Lusztig nilpotent variety}}
Let $T^*\mathcal{M}_{Q}=\bigsqcup_{\d}T^*\mathcal{M}_{Q,\d}$ be the cotangent\footnote{to be precise, we only consider the underived, or $H^0$-truncation, of the cotangent stacks here} stack of $\mathcal{M}_{Q}$. This may be realized explicitly as follows. Let $\overline{Q}=(I, \Omega \sqcup \Omega^*)$ be the double of $Q$, obtained by replacing every arrow $h$ in $Q$ by a pair $(h,h^*)$ of arrows going in opposite directions. The preprojective algebra $\Pi_Q$ is the quotient of the path algebra $\k\overline{Q}$ by the two-sided ideal generated by $\sum_{h \in \Omega} [h, h^*]$. Unless $Q$ is of finite type, the abelian category $\text{Rep}_{\k}\Pi_Q$ is of global dimension two. The stack of $\k$-representations of $\Pi_Q$ is naturally identified with $T^*\mathcal{M}_Q$. We say that a representation $M$ of $\Pi_Q$ is \textit{nilpotent} if there exists a filtration $M \supset M_1 \supset \cdots \supset M_l=\{0\}$ for which $\Pi^+_Q (M_i) \subseteq M_{i+1}$, where $\Pi^+_Q \subset \Pi_Q$ is the augmentation ideal. Following Lusztig, we define the \textit{nilpotent variety} (or stack) $\Lambda_Q = \bigsqcup_{\d} \Lambda_{Q,\d} \subset T^*\mathcal{M}_Q$ as the substack of nilpotent representations of $\Pi_Q$.

\begin{theorem}[Kashiwara-Saito, \cite{KaSa}] The stack $\Lambda_Q$ is a lagrangian substack of $T^*\mathcal{M}_Q$ and
%and for any perverse sheaf $\mathbb{P} \in \mathcal{K}_Q$ we have $SS(\mathbb{P}) \subset \Lambda_Q$. Moreover, 
for any $\d \in \N^I$ we have
$\# Irr(\Lambda_{Q,\d})=\dim\; U(\mathfrak{n}_Q)[\d]$.
\end{theorem}

The above theorem is strongly motivated by Lusztig's geometric lift of the spherical Hall algebra to a category of perverse sheaves on $\mathcal{M}_Q$, yielding the \textit{canonical basis} $\mathbf{B}$ of $U_v(\mathfrak{n}_Q)$, see \cite{LusJAMS}. It turns out that $\Lambda_Q$ is precisely the union of the singular support of Lusztig's simple perverse sheaves, hence the above may be seen as a microlocal avatar of Lusztig's construction. 
Kashiwara and Saito prove much more, namely they equip $Irr(\Lambda)$ with the combinatorial structure of a \textit{Kashiwara crystal} which they identify as the crystal $B(\infty)$ of $U(\mathfrak{n}_Q)$. This also yields a canonical bijection between $\mathbf{B}$ and $Irr(\Lambda_Q)$.

\vspace{.1in}

\paragraph{\textbf{2.5.\;Arbitrary quivers and nilpotent Kac polynomials}} What happens when the quiver $Q$ \textit{does} have edge loops, such as the $g$-loop quiver $S_g$ ? This means that the matrix of the Euler form may have some (even) nonpositive entries on the diagonal, hence it is associated to a \textit{Borcherds} algebra rather than to a Kac-Moody algebra. Accordingly, we call \textit{real}, resp. \textit{isotropic}, resp. \textit{hyperbolic} a vertex carrying zero, resp. one, resp. at least two edge loops, and in the last two cases we say that the vertex is \textit{imaginary}. One may try to get Hall-theoretic constructions of the Borcherds algebra $\g_Q$ associated to $Q$ (see e.g. \cite{KangSch}, \cite{KangKashSch}) but the best thing to do in order to get a picture as close as possible to the one in \textbf{2.1}--\textbf{2.4.} seems to be to consider instead a slightly larger algebra $\g_Q^B$ defined by Bozec which has as building blocks the usual $\mathfrak{sl}_2$ for real roots, the Heisenberg algebra $H$ for isotropic roots and a free Lie algebra $H'$ with one generator in each degree for hyperbolic roots, see \cite{Bozec1}. 

The Hall algebra $\widetilde{\H}_Q$ is defined just as before, but we now let $\widetilde{\H}^{\text{sph}}_Q$ be the subalgebra generated
by elements $1_{S_i}$ for $i$ real vertices and by the characteristic functions $1_{\mathcal{M}_{l\epsilon_i}}$ for $i$ imaginary and $l \in \N$. The stack of nilpotent representations of $\Pi_Q$ is not lagrangian anymore in general, we consider instead the stack of \textit{seminilpotent} representations, i.e. representations $M$ for which there exists a filtration $M \supset M_1 \supset \cdots \supset M_l=\{0\}$ such that $h(M_i / M_{i+1})=0$ for $h \in \Omega$ or $h \in \Omega^*$ \textit{not} an edge loop.

\begin{theorem}[Bozec, \cite{Bozec1},\cite{Bozec2}, Kang, \cite{KangBozec2}]\label{T:B1} The following hold~:
\begin{enumerate}
\item[i)] $\widetilde{\H}^{\text{sph}}_Q$ is isomorphic to $U_v(\mathfrak{b}^B_Q)$,
%\item[ii)] $K_0^{gr}(\mathcal{K}_Q)$ is isomorphic to $U_v(\mathfrak{n}^B_Q)$,
\item[ii)] The stack $\Lambda_Q$ of seminilpotent representations is lagrangian in $T^*\mathcal{M}_Q$. For any $\d \in \N^I$ we have $\#Irr(\Lambda_{Q,\d})=\dim\;U(\mathfrak{n}_Q^B)[\d]$.
\end{enumerate}
\end{theorem}

The theories of canonical bases and crystal graphs also have a natural extension to this setting. What about the relation to Kac polynomials ? Here as well we need some slight variation; let us call \textit{$1$-nilpotent} a representation $M$ of $Q$ in which the edge loops at any given imaginary vertex $i$ generate a nilpotent (associative) subalgebra of $\text{End}(M_i)$. Of course any representation of a quiver with no edge loops is $1$-nilpotent.

\begin{theorem}[Bozec-S.-Vasserot, \cite{BSV}]\label{T:BSV} For any $\d$ there exists a (unique) polynomial $A^{\text{nil}}_{Q,\d} \in \mathbb{Z}[t]$, independent of the orientation of $Q$,
such that for any finite field $\k$ the number of isomorphism classes of absolutely indecomposable $1$-nilpotent $\k$-representations of $Q$ of dimension $\d$ is equal to $A^{\text{nil}}_{Q,\d}(\#\k)$. Moreover, for any $\d \in \N^I$ we have
\begin{enumerate}
\item[i)] $A^{\text{nil}}_{Q,\d}(1)=A_{Q,\d}(1)$,
\item[ii)]$A^{\text{nil}}_{Q,\d}(0)=\dim\;\mathfrak{g}^B_Q[\d]$.
\end{enumerate}
\end{theorem}

We will sketch a proof of the fact that $A^{\text{nil}}_{Q,\d}(t) \in \N[t]$ in Section~\textbf{4.4}.

\medskip

\section{Quivers, Kac polynomials and graded Borcherds algebras I \\--counting cuspidals--}

\medskip

Starting form a Kac-Moody algebra $\g$ (or its variant defined by Bozec) we build a quiver $Q$ by orienting the Dynkin diagram of $\g$ in any fashion, and obtain a realization of $U_v(\mathfrak{b})$ as the spherical Hall algebra $\widetilde{\H}_Q^{\text{sph}}$; now, we may ask the following question~: what is the relation between $\g$ and the \textit{full} Hall algebra ? Here's a variant of this question~: Kac's constant term conjecture gives an interpretation of the constant term of the Kac polynomial $A_{Q,\d}(t)$; is there a similar Lie-theoretic meaning for the \textit{other} coefficients~?
As we will explain in \textbf{3.3.} below, the two questions turn out to have a beautiful common conjectural answer~ : the full Hall algebra $\widetilde{\H}_Q$ is related to a \textit{graded} Borcherds algebra $\widetilde{\g}_Q$ whose \textit{graded} multiplicities $\dim_{\Z}\;\widetilde{\g}_Q[\d]=\sum_l \dim\;\widetilde{\g}_Q[\d,l]t^l$ are equal to $A_{Q,\d}(t)$. What's more, $\widetilde{\g}_Q$, like $\g_Q$, is independent of the choice of the orientation of the quiver and is hence \textit{canonically} attached to $\g_Q$. Although the above makes sense for an arbitrary quiver $Q$, there is an entirely similar story in the nilpotent setting (better suited for a quiver with edge loops), replacing $A_{Q,\d}(t)$ by $A^{\text{nil}}_{Q,\d}(t)$ and $\H_Q$ by the subalgebra $\H^{\text{nil}}_Q$ of functions on the stack $\mathcal{M}^{\text{nil}}_{Q}$ of $1$-nilpotent representations of $Q$. This is then conjecturally related to a graded Borcherds algebra $\widetilde{\mathfrak{g}}_Q^{\text{nil}}$.

\vspace{.1in}

\paragraph{\textbf{3.1. The full Hall algebra.}} The first general result concerning the structure of the full Hall algebra $\widetilde{\H}_Q$ is due to Sevenhant and Van den Bergh. Let us call an element $f \in \H_Q$ is \textit{cuspidal} if $\Delta(f) =f \otimes 1 + k_f \otimes 1$. Denote by $\H_Q^{\text{cusp}}=\bigoplus_{\d}\H_Q^{\text{cusp}}[\d]$ the space of cuspidals, and set
$C_{Q/\k,\d}=\dim\;\H_Q^{\text{cusp}}[\d]$. We write $Q/\k$ instead of $Q$ to emphasize the dependance on the field. Consider the infinite Borcherds-Cartan data $(A_{Q/\k}, \mathbf{m}_{Q/\k})$, with $\N^I \times \N^I$ Cartan matrix $A_{Q/\k}=(a_{\d,\d'})_{\d,\d'}$ and charge function $\mathbf{m}_{Q/\k}:~ \N^I \to \N$ defined as follows~:
$$a_{\d,\d'}=(\d,\d'), \qquad \mathbf{m}_{Q/\k}(\d)=C_{Q/\k,\d} \qquad \forall\; \d,\d' \in \N^I.$$
Denote by $\widetilde{\g}_{Q/\k}$ Borcherds algebra associated to $(A_{Q/\k},\mathbf{m}_{Q/\k})$.

\begin{theorem}[Sevenhant-Van den Bergh, \cite{SVdB}] The Hall algebra ${\H}_{Q/\k}$ is isomorphic to the positive nilpotent subalgebra $U_v(\widetilde{\mathfrak{n}}_{Q/\k})$ of the Drinfeld-Jimbo quantum enveloping algebra of $\widetilde{\g}_{Q/\k}$.
\end{theorem}

In other words, the cuspidal element span the spaces of \textit{simple root vectors} for the Hall algebra $\H_{Q/\k}$. Note that adding the 'Cartan' $\mathbf{K}$ to $\H_Q$ will only produce a quotient of $U_v(\widetilde{\mathfrak{b}}_{Q/\k})$, whose own Cartan subalgebra is not finitely generated. 

\vspace{.1in}

\paragraph{\textbf{3.2.\;Counting cuspidals for quivers.}} To understand the structure of $\widetilde{\g}_{Q/\k}$ better, we first need to understand the numbers $C_{Q/\k,\d}$. In this direction, we have the following theorem, proved by Deng and Xiao for quivers with no edge loops, and then extended to the general case in \cite{BScusp}~:

\begin{theorem}[Deng-Xiao, \cite{DX}, Bozec-S., \cite{BScusp}] For any $Q=(I,\Omega)$ and $\d \in \N^I$ there exists a (unique) polynomial $C_{Q,\d}(t) \in \mathbb{Q}[t]$ such that for any finite field $\k$, $C_{Q/\k,\d}=C_{Q,\d}(\#\k)$.
\end{theorem}

Let us say a word about the proof. By the Krull-Schmidt theorem, there is an isomorphism of $\N^I$-graded vector spaces
between $\H_{Q/\k}$ and $Sym(\bigoplus_{\d} Fun(\mathcal{M}_{Q/\k,\d}^{ind}(\k),\C))$, where $\mathcal{M}_{Q/\k,\d}^{ind}(\k)$
is the groupoid of indecomposable $\k$-representations of $Q$ of dimension $\d$. This translates into the equality of generating series
\begin{equation}\label{E:cuspgenseries}
\sum_{\d} (\dim\;\H_{Q/\k}[\d])z^d = \text{Exp}_z\left( \sum_{\d} I_{Q,\d}(t)z^d\right)_{|t=\#\k}
\end{equation}
where $\text{Exp}_z$ stands for the plethystic exponential with respect to the variable $z$, see e.g. \cite{Mozpleth}.
From the PBW theorem, 
\begin{equation}\label{E:cuspgenseries2}
\sum_{\d} (\dim\;U(\widetilde{\mathfrak{n}}_{Q/\k})[\d])z^d = \text{Exp}_z\left( \sum_{\d} (\dim\;\widetilde{\mathfrak{n}}_{Q/\k}[\d])z^d\right)
\end{equation}
from which it follows that $\dim\;\widetilde{\mathfrak{n}}_{Q/\k}[\d]=I_{Q,\d}(\#\k)$ for any $\d$. We are in the following situation~: we know the character\footnote{there is a subtle point here : $\widetilde{\mathfrak{g}}_{Q/\k}$ is $\Z^{(\N^I)}$-graded, but we only know the character after projection of the grading to $\N^I$. } of the Borcherds algebra $\widetilde{\mathfrak{g}}_{Q/\k}$ and from that we want to determine its Borcherds-Cartan data. This amounts to inverting the Borcherds character (or denominator) formula, which may be achieved in each degree by some finite iterative process\footnote{here we use the fact that the Weyl group of $\widetilde{\mathfrak{g}}_{Q/\k}$ only depends on the subset of real simple roots, which is the same as that of $\g_Q$.}. Following through this iterative process and using the fact that $I_{Q,d}(\#\k)$ is a polynomial in $\#\k$, one checks that the Borcherds-Cartan data, in particular the charge function $\mathbf{m}_{Q/\k}(\d)=C_{Q/\k,\d}$ is also polynomial in $\#\k$. This process is constructive, but as far as we know no closed formula is known in general (see \cite{BScusp} for some special cases, such as totally negative quivers like $S_g, g >1$). And in particular, the structure (and the size) of $\H_{Q/\k}$ depends heavily on $\k$.

\vspace{.1in}

\paragraph{\textbf{3.3.\;Counting absolutely cuspidals for quivers}} In order to get a more canonical structure out of $\H_{Q/\k}$ one is tempted to view the polynomial $C_{Q,\d}(t)$ as giving the graded dimension of some vector space. Unfortunately, as can readily be seen for the Jordan quiver $S_1$ and $\d=2$, the polynomial $C_{Q,\d}(t)$ fails to be integral or positive in general. This is familiar~: the same things happens for Kac's $A$ and $I$ polynomials and one might guess that the better thing to count would be \textit{absolutely} cuspidal elements of $\H_{Q/\k}$. Unfortunately (and contrary to the case of curves), there is at the moment no known definition of an absolutely cuspidal element of $\H_{Q/\k}$ ! One way around this is to use a well-known identity between Kac $A$-polynomials and $I$-polynomials to rewrite \eqref{E:cuspgenseries}
as 
\begin{equation}\label{E:cuspgenseries3}
\sum_{\d} \dim(U(\widetilde{\mathfrak{n}}_{Q/\k})[\d])z^d=\sum_{\d} \dim(H_{Q/\k}[\d])z^d = \text{Exp}_{z,t}\left( \sum_{\d} A_{Q,\d}(t)z^d\right)_{|t=\#\k}
\end{equation}
and interpret the right hand side as the character of the enveloping algebra of (the nilpotent subalgebra of) a putative \textit{graded} Borcherds algebra $\widetilde{\g}_Q$ whose graded character is given by the polynomials $A_{Q,\d}(t)$
\begin{equation}\label{E:gradedimBKM}
\text{dim}_{\Z}\; \widetilde{\mathfrak{g}}_Q[\d]=A_{Q,\d}(t).
\end{equation}
Taking this as a working hypothesis, we can run through the above iterative process (in the graded sense) to get a well-defined family of polynomials
$C^{\text{abs}}_{Q,\d}(t)$.

\begin{theorem}[Bozec-S., \cite{BScusp}] For any $Q$ and $\d$ we have $C^{\text{abs}}_{Q,\d}(t) \in \Z[t]$. In addition,  we have
$C_{Q,\d}^{\text{abs}}(t)=C_{Q,\d}(t)$ for any hyperbolic $\d$ (i.e. $(\d,\d)<0$) while for any indivisible isotropic $\d$,
$$\textup{Exp}_z\left(\sum_{l \geq 1} C_{Q,l\mathbf{d}}(t)z^{l}\right)=\textup{Exp}_{t,z}\left(\sum_{l \geq 1} C^{\text{abs}}_{Q,l\mathbf{d}}(t)z^{l}\right)$$
\end{theorem}

Our definition of $C_{Q,\d}^{\text{abs}}(t)$ was motivated by the putative existence of a graded Borcherds algebra $\widetilde{\g}_Q$. This existence is \textit{a posteriori} equivalent to the following conjecture~:

\begin{conjecture} For any $Q$ and $\d$ we have $C_{Q,\d}^{\text{abs}}\in \N[t]$.
\end{conjecture}

We stress that the above positivity conjecture is strictly stronger than the positivity of Kac's $A$-polynomials. Indeed by construction Kac polynomials are integral positive polynomials in the $C_{Q,\d}^{\text{abs}}(t)$. We will give in Sections~\textbf{4.4, 4.6.} two other (conjectural) constructions of the Lie algebra $\widetilde{\g}_Q$, based on the geometry of $T^*\mathcal{M}_Q, \Lambda_Q$ and on Nakajima quiver varieties respectively, either of which would imply the above positivity conjecture. We finish this section with another (somewhat imprecise) conjecture, which is a quiver analog of a conjecture of Kontsevich and Deligne, see Section~\textbf{6.1}.

\begin{conjecture}\label{C:geomcusp} For any $Q,\d$ there exists a 'natural' algebraic variety $\mathcal{C}_{Q,\d}$ defined over $\mathbb{Z}$ such that for any finite field $\k$ we have $C_{Q,\d}^{\text{abs}}(\#\k)=\#\mathcal{C}_{Q,\d}(\k)$.
\end{conjecture}
There are obvious variants of this conjecture, replacing $\mathcal{C}_{Q,\d}$ by a complex algebraic variety and the point count by the Poincar\'e polynomial, fixing the characteristic of $\k$, etc.
 
\smallskip
 
\noindent
\textit{Examples~:} For $Q$ an affine quiver, we have $C_{Q,\d}^{\text{abs}}(t)=1,t$ or $0$ according to whether $\d \in \{\epsilon_i\}_i$, $\d \in \N\delta$ for $\delta$ the indivisible imaginary root, or neither of the above. For $Q=S_3$ and $\d=3$ we have
$$C^{\text{abs}}_{S_g,3}=\frac{t^{9g-3}-t^{5g+2}-t^{5g-2}-t^{5g-3}+t^{3g+2}+t^{3g-2}}{(t^2-1)(t^3-1)}.$$
 
\vspace{.1in}

\paragraph{\textbf{3.4.\;More questions than answers}} As we have seen, there is a method for counting cuspidal (or even 'absolutely' cuspidal) functions, via Kac polynomials and the Borcherds character formula. But is it possible to get a geometric parametrization of these cuspidals (either in the sense of conjecture~\ref{C:geomcusp}, or as in the Langlands program in terms of some 'spectral data' \footnote{a preliminary question~: what plays the role of Hecke operators in this context?}) ? Is it possible to explicitly construct cuspidal functions ? Can one lift a suitable basis of the space of cuspidal functions to some perverse sheaves on $\mathcal{M}_Q$ ? In other words, is there a theory of canonical bases for $U_v(\widetilde{\n}_Q)$ ? If so, can one describe explicitly the Ext-algebra of that category of perverse sheaves~? Equivalently, is there a analog, for $\widetilde{\g}_Q$, of the Khovanov-Lauda-Rouquier algebra, see \cite{KLR1}, \cite{KLR2} ?

 \vspace{.1in}

\section{Quivers, Kac polynomials and graded Borcherds algebras II \\
--Cohomological Hall algebras and Yangians--}
%2,5 pages
\medskip

The previous section offered a conjectural definition of Lie algebras $\widetilde{\g}_Q, \widetilde{\g}^{\text{nil}}_Q$ whose graded multiplicities are $A_{Q,\d}(t), A_{Q,\d}^{\text{nil}}(t)$, but the construction --involving something like a 'generic form' for the full Hall algebras $\H_Q, \H^{\text{nil}}_Q$-- was somewhat roundabout. In this section we describe a geometric construction, in terms of the cohomology of $T^*\mathcal{M}_Q$ or $\Lambda_Q$, of algebras --the \textit{cohomological Hall algebras}-- which are deformations of the envelopping algebras $U(\bar{\mathfrak{n}}_Q[u]), U(\bar{\mathfrak{n}}^{\text{nil}}_Q[u])$ for certain graded Lie algebras $\bar{\mathfrak{n}}_Q, \bar{\mathfrak{n}}^{\text{nil}}_Q$ satisfying $\dim_{\Z}\;\bar{\mathfrak{n}}_Q=A_{Q,\d}(t), \dim_{\Z}\;\bar{\mathfrak{n}}^{\text{nil}}_Q=A^{\text{nil}}_{Q,\d}(t)$.
Of course, it is expected that $\bar{\mathfrak{n}}_Q, \bar{\mathfrak{n}}^{\text{nil}}_Q$ are positive halves of Borcherds algebras $\bar{\mathfrak{\g}}_Q, \bar{\mathfrak{\g}}^{\text{nil}}_Q$ and thus coincide with $\widetilde{\mathfrak{g}}_Q, \widetilde{\mathfrak{g}}^{\text{nil}}_Q$. One nice output of this construction is that it yields for free a whole family of representations (in the cohomology of Nakajima quiver varieties). In this section, $\k=\C$.

\vspace{.1in}

\paragraph{\textbf{4.1.\;Geometry of the nilpotent variety}}  Let us begin by mentioning a few remarkable geometric properties of the (generally very singular) stacks $T^*\mathcal{M}_Q$, $\Lambda_{Q}$ introduced in ~\textbf{2.4.} and \textbf{2.5.}. Here $Q=(I,\Omega)$ is an arbitrary quiver. Denote by $T=(\C^*)^2$ be the two-dimensional torus acting on $\Lambda_Q$ and $T^*\mathcal{M}_Q$ by scaling the arrows in $\Omega$ and $\Omega^*$ respectively. 

\begin{theorem}[S.-Vasserot, \cite{SVCOHA1}, Davison, \cite{Davison}] The $T$-equivariant Borel-Moore homology space $H_*^T(\Lambda_Q,\Q)=\bigoplus_{\d} H_*^T(\Lambda_{Q,\d}, \Q)$ is even, pure, and free over $H^*_T(pt, \Q)$.  The same holds for $T^*\mathcal{M}_Q$.
\end{theorem}

We will only consider (co)homology with $\Q$-coefficients and drop the $\Q$ in the notation. The above theorem for $\Lambda_Q$ is proved by constructing a suitable compactification of $\Lambda_{Q,\d}$ for each $\d$, itself defined in terms of Nakajima varieties, see \cite{SVCOHA1}, \cite{Nak}. The second part of the theorem relies on dimensional reduction from a $3$d
Calabi-Yau category, see \cite{Davison}. Kac polynomials are intimately related to the cohomology of $\Lambda_Q$ and $T^*\mathcal{M}_Q$. More precisely,

\begin{theorem}[Bozec-S.-Vasserot, Davison, Mozgovoy]\label{T:BSV2}The Poincar\'e polynomials of $\Lambda_Q$ and $T^*\mathcal{M}_Q$ are respectively given by~:
\begin{equation}\label{E:PP1}
\sum_{i,\d} \dim\; H_{2i}(\Lambda_{Q,\d})t^{\langle \d,\d\rangle+i}z^\d= \text{Exp}\left(\sum_{\d} \frac{A^{\text{nil}}_{Q,\d}(t^{-1})}{1-t^{-1}}z^{\d}\right)
\end{equation}
\begin{equation}\label{E:PP2}
\sum_{i,\d} \dim\; H_{2i}(T^*\mathcal{M}_{Q,\d})t^{\langle \d,\d\rangle+i}z^\d= \text{Exp}\left(\sum_{\d} \frac{A_{Q,\d}(t)}{1-t^{-1}}z^{\d}\right).
\end{equation}
\end{theorem}
Equality (\ref{E:PP1}) is proved in \cite{BSV} and hinges again on some partial compactification of the stacks $\Lambda_{Q,\d}$ defined in terms of Nakajima quiver varieties. The second equality (\ref{E:PP2}) is obtained by combining the purity results of \cite{Davison} with the point count of \cite{Mozpleth}. In some sense, (\ref{E:PP1}) and (\ref{E:PP2}) are Poincar\'e dual to each other. Note that
$\dim\; \Lambda_{Q,\d}=-\langle \d,\d\rangle$ hence (\ref{E:PP1}) is a series in $\C[[t^{-1}]][z]$. Taking the constant term in $t^{-1}$ yields
$\sum_{\d} \#Irr(\Lambda_{\d})z^{\d}=Exp(\sum_{\d} A_{\d}^{\text{nil}}(0)z^\d)$, in accordance with Theorem~\ref{T:B1} ii), Theorem~\ref{T:BSV}
ii). In other terms, (\ref{E:PP1}) provides a geometric interpretation (as the Poincar\'e polynomial of some stack) of the \textit{full}\; Kac polynomial $A^{\text{nil}}_{Q,\d}(t)$ rather than just its constant term and passing from $A^{\text{nil}}_{Q,\d}(0)$ to $A_{Q,\d}^{\text{nil}}(t)$ essentially amounts to passing from $H_{top}(\Lambda_{Q,\d})$ to $H_*(\Lambda_{Q,\d})$. Of course, (\ref{E:PP2}) has a similar interpretation.

\vspace{.1in}

\paragraph{\textbf{4.2.\;Cohomological Hall algebras}} The main construction of this section is the following~:

\begin{theorem}[S.-Vasserot, \cite{SVIHES}, \cite{SVCOHA1}] The spaces $H_*^T(T^*\mathcal{M}_Q)$ and
$H_*^T(\Lambda_{Q})$ carry natural $\mathbb{Z} \times \N^I$-graded algebra structures. The direct image morphism $i_*: H^T_*(\Lambda_Q) \to H^T_*(T^*\mathcal{M}_Q)$ is an algebra homomorphism.
\end{theorem}

Let us say a few words about how the multiplication map $H^T_*(T^*\mathcal{M}_{Q,\d_1}) \otimes H^T_*(T^*\mathcal{M}_{Q,\d_2}) \to H^T_*(T^*\mathcal{M}_{Q,\d_1+\d_2})$ is defined. There is a convolution diagram
$$\xymatrix{T^*\mathcal{M}_{Q,\d_1} \times T^*\mathcal{M}_{Q,\d_2} & Z_{\d_1,\d_2} \ar[l]_-{q} \ar[r]^-{p} & T^*\mathcal{M}_{Q,\d_1+\d_2}}$$
similar to (\ref{E:conv}), where $Z_{\d_1,\d_2}$ is the stack of inclusions $\bar{M} \subset \bar{R}$ with $\bar{M}, \bar{R}$ representations of $\Pi_Q$ of respective dimensions $\d_2, \d_1+\d_2$, and $p$ and $q$ are the same as in (\ref{E:conv}). The map $p$ is still proper, so that $p_*~: H^T_*(Z_{\d_1,\d_2}) \to H^T_*(T^*\mathcal{M}_{Q,\d_1+\d_2})$ is well-defined, but $q$ is not regular anymore and we cannot directly define a Gysin map $q^*$. Instead, we embedd $Z_{\d_1,\d_2}$ and $T^*\mathcal{M}_{Q,\d_1} \times T^*\mathcal{M}_{Q,\d_2}$ into suitable \textit{smooth} moduli stacks of representations of the path algebra of the double quiver $\overline{Q}$ (without preprojective relations)
and define a refined Gysin map $q^!~:H^T_*(T^*\mathcal{M}_{Q,\d_1} \times T^*\mathcal{M}_{Q,\d_2}) \to H^T_*(Z_{\d_1,\d_2})$. The multiplication map is then $m=p_*q^!$. Note that it is of cohomological degree $-\langle \d_1, \d_2\rangle - \langle \d_2,\d_1\rangle$; to remedy this, we may (and will) shift the degree of $H_*^T(T^*\mathcal{M}_{Q,\d})$ by $\langle\d,\d\rangle$.
There is another notion of cohomological Hall algebra due to Kontsevich and Soibelman, associated to any Calabi-Yau category of dimension \textit{three} (\cite{KonSoi}); as shown by Davison \cite{DavisonKS=SV} (see also Yang-Zhao \cite{YY2}) using a dimensional reduction argument, the algebras $H^T_*(T^*\mathcal{M}_Q), H^T_*(\Lambda_Q)$ arise in that context as well.  In addition, the above construction of the cohomological Hall algebras can be transposed to any oriented cohomology theory which has proper pushforwards and refined Gysin maps, such as Chow groups, K-theory, elliptic cohomology or Morava K-theory (see \cite{YZ}). This yields in a uniform way numerous types of quantum algebras with quite different flavors. 

\smallskip

Not very much is known about the precise structure of $H^T_*(T^*\mathcal{M}_Q)$ or $H^T_*(\Lambda_Q)$ in general (and even less so for other types of cohomology theories).  However, it is possible to give a simple set of generators; for any $\d$ there is an embedding of
$\mathcal{M}_{Q,\d}$ in $\Lambda_{Q,\d}$ as the zero-section, the image being one irreducible component of $\Lambda_{Q,\d}$. 

\begin{theorem}[S.-Vasserot, \cite{SVCOHA1}]\label{T:generation} The algebra $H^T_*(\Lambda_Q)$ is generated over $H^*_T(pt)$ by the collection of subspaces
$$\{H^T_*(\mathcal{M}_{Q,\epsilon_i})\;|\; i \in I^{re} \cup I^{iso} \}\cup \{H^T_*(\mathcal{M}_{Q,l\epsilon_i})\;|\; i \in I^{hyp}, l \in \N \}$$
\end{theorem}
Note that $\mathcal{M}_{Q,\d} = E_{\d}/G_{\d}$, where $E_{\d}$ is a certain representation of $G_{\d}=\prod_i GL(\d_i)$, hence 
$H^T_*(\mathcal{M}_{Q,\d})=H^*_{G_{\d} \times T}(pt) \cdot [\mathcal{M}_{Q,\d}] \simeq \Q[q_1,q_2,c_{l}(M_i)\;|\; i \in I, l \leq \d_i]$, where $q_1,q_2$ are the equivariant parameters corresponding to $T$. The above result is for the nilpotent stack $\Lambda_{Q,\d}$, but one can show that $H^T_*(\Lambda_{Q}) \otimes \Q(q_1,q_2) =H^T_*(T^*\mathcal{M}_Q)\otimes \Q(q_1,q_2)$ \footnote{in other words $H^T_*(\Lambda_Q)$ and $H^T_*(T^*\mathcal{M}_Q)$ are two different integral forms of the same $\Q(q_1,q_2)$-algebra; this explains the discrepancy between usual and nilpotent Kac polynomials.} so that the same generation result holds for $H^T_*(T^*\mathcal{M}_Q)\otimes \Q(q_1,q_2)$. 

\vspace{.1in}

\paragraph{\textbf{4.3. Shuffle algebras.}} One can also give an algebraic model for a certain \textit{localized} form of $H^T_*(T^*\mathcal{M}_Q)$ (or $H^T_*(\Lambda_Q)$). Namely, let $\overline{E}_{\d}=T^*E_{\d}$ be the vector space of representations of $\overline{Q}$ in $\bigoplus_{i} \k^{\d_i}$. The direct image morphism $i_*: H^{T}_*(T^*\mathcal{M}_{Q,\d}) \to H^T_*(\overline{E}_{\d}/G_{\d})$ is an isomorphism after tensoring by $Frac(H^*_{T \times G_{\d}}(pt))$. What's more, it is possible to equip 
$$Sh^{H^T_*}_Q:=\bigoplus_{\d}H^T_*(\overline{E}_{\d}/G_{\d})\simeq \bigoplus_{\d} H^*_{T \times G_{\d}}(pt)$$ with the structure of an associative algebra (such that $i_*$ becomes an algebra morphism), described explicitly as a \textit{shuffle algebra}. More precisely, let us identify 
\begin{equation}\label{E:shuffle10}
H^*_{T \times G_\d}(pt) \simeq \Q[q_1,q_2, z_{i,l}\;|\; i \in I,  l \leq \d_i]^{W_{\d}}
\end{equation}
where $W_{\d}=\prod_i \mathfrak{S}_{\d_i}$ is the Weyl group of $G_{\d}$. To unburden the notation, we will collectively denote the variables $z_{i,1}, \ldots, z_{i,\d_i}$ (for all $i \in I$) by $\underline{z}_{[1,\d]}$, and use obvious variants of that notation. Also, we will regard $q_1,q_2$ as scalars and omit them from the notation. Fix dimension vectors $\d,\e$ and put $\n=\d+\e$. For two integers $r,s$ we denote by $Sh_{r,s} \subset \mathfrak{S}_{r+s}$ the set of $(r,s)$-shuffles, i.e. permutations $\sigma$ satisfying $\sigma(i) < \sigma(j)$ for $1 \leq i <j \leq r$ and $r< i<j \leq r+s$, and we put
$Sh_{\d,\e}=\prod_i Sh_{\d_i,\e_i} \subset \prod_i \mathfrak{S}_{\n_i}=W_{\n}$. In terms of (\ref{E:shuffle10}), the multiplication map $H^*_{T \times G_{\d}}(pt) \otimes H^*_{T \times G_{\e}}(pt) \to H^*_{T \times G_{\n}}(pt)$ now reads
$$(f * g)(\underline{z}_{[1,\n]})=\sum_{\sigma \in Sh_{\d,\e}}\sigma \left[K_{\d,\e}(\underline{z}_{[1,\n]})\cdot f(\underline{z}_{[1,\d]}) \cdot g(\underline{z}_{[\d+1,\n]})\right]$$
where $K_{\d,\e}(\underline{z}_{[1,\n]})=\prod_{s=0}^2K^{(s)}_{\d,\e}(\underline{z}_{[1,\n]})$ with
$$K^{(0)}_{\d,\e}(\underline{z}_{[1,\n]})=\prod_{i \in I}\hspace{-.05in}\prod_{\substack{1 \leq l \leq \d_i\\ \d_i+1 \leq k  \leq \n_i}} \hspace{-.15in}(z_{i,l}-z_{i,k})^{-1},$$
$$K^{(1)}_{\d,\e}(\underline{z}_{[1,\n]})=\prod_{h \in \Omega}\prod_{\substack{1 \leq l \leq \d_{h'}\\ \d_{h''}+1 \leq k  \leq \n_{h''}}} \hspace{-.15in}(z_{h'\hspace{-0.03in},l}-z_{h''\hspace{-0,03in},k}-q_1) 
\prod_{\substack{1 \leq l \leq \d_{h''}\\ \d_{h'}+1 \leq k  \leq \n_{h'}}} \hspace{-.15in}(z_{h''\hspace{-0,03in},l}-z_{h'\hspace{-0,03in},k}-q_2) $$
and
$$K^{(2)}_{\d,\e}(\underline{z}_{[1,\n]})=\prod_{i \in I}\hspace{-.05in}\prod_{\substack{1 \leq l \leq \d_i\\ \d_i+1 \leq k  \leq \n_i}} \hspace{-.15in}(z_{i,k}-z_{i,l}-q_1-q_2).$$

As shown in \cite{SVCOHA1}, $H^T_*(T^*\mathcal{M}_{Q,\d})$ and $H^T_{*}(\Lambda_{Q,\d})$ are torsion-free and of generic rank one as modules over $H^*_{T \times G_\d}(pt)$, hence the localization map $i_*$ is injective and Theorem~\ref{T:generation} yields a description of the cohomological Hall algebra $H_*^T(\Lambda_Q)$ as a subalgebra of the above shuffle algebra, generated by an explicit collection of polynomials. This allows one to identify the rational form of $H_*^T(\mathcal{M}_Q)$ with the positive half of the Drinfeld Yangian $Y_h(\mathfrak{g}_Q)$ when $Q$ is of finite type and with the positive half of the Yangian version of the elliptic Lie algebra $\mathfrak{g}_{Q_0}[s^{\pm 1}, t^{\pm 1}] \oplus K$ when $Q$ is an affine Dynkin diagram. Here $K$ is the full central extension of the double loop algebra $\mathfrak{g}_{Q_0}[s^{\pm 1}, t^{\pm 1}] $.
Beyond these case, shuffle algebras tend to be rather difficult to study and the algebraic structure of $H_*^T(\mathcal{M}_Q)$ (or $H_*^T(\Lambda_Q)$) is still mysterious (see, however \cite{Negut1}, \cite{NegutAGT} for 
some important applications to the geometry of instanton moduli spaces in the case of the Jordan or affine type $A$ quivers). 

\vspace{.1in}

There are analogous shuffle algebra models in the case of an arbitrary oriented Borel-Moore homology theory, but the torsion-freeness statement remains conjectural in general. In order to give the reader some idea of what these are, as well as for later use, let us describe the localized K-theoretic Hall algebra of the $g$-loop quiver $S_g$. In this case, there is a $g+1$-dimensional torus $T_g$ acting by rescaling the arrows :
$$(\xi_1, \xi_2, \ldots, \xi_g, p) \cdot (h_1, h_1^*,\ldots, h_g, h_g^*) = (\xi_1h_1, p\xi_1^{-1}h_1^*, \ldots, \xi_g h_g, p\xi_g^{-1} h_g^*)$$
and we have an identification
$$Sh^{K^{T_g}}_{S_g}:=\bigoplus_{d} K^{T_g \times G_d}(pt)\simeq \bigoplus_{d}\mathbb{Q}[\xi_1^{\pm 1}, \ldots, \xi^{\pm 1}_g,p^{\pm 1}; z_1^{\pm 1}, \ldots, z_d^{\pm 1}]^{\mathfrak{S}_d}.$$
The multiplication takes the form
\begin{equation}\label{E:shuffle1}
(f * g)(\underline{z}_{[1,n]})=\sum_{\sigma \in Sh_{d,e}}\sigma \left[K_{d,e}(\underline{z}_{[1,n]})\cdot f(\underline{z}_{[1,d]}) \cdot g(\underline{z}_{[d+1,n]})\right]
\end{equation}
where $K_{d,e}(\underline{z}_{[1,n]})=\prod_{s=0}^2K^{(s)}_{d,e}(\underline{z}_{[1,n]})$ with
$$K_{d,e}^{(0)}(\underline{z}_{[1,n]})=\hspace{-.1in}\prod_{\substack{1 \leq l \leq d \\ d+1 \leq k \leq n}} (1-z_l/z_k)^{-1},\qquad K_{d,e}^{(1)}(\underline{z}_{[1,n]})=\prod_{\substack{1 \leq l \leq d \\ d+1 \leq k \leq n}} \prod_{u=1}^g(1-\xi_u^{-1}z_l/z_k)(1-p^{-1}\xi_uz_l/z_k)$$
$$K_{d,e}^{(2)}(\underline{z}_{[1,n]})=\hspace{-.1in}\prod_{\substack{1 \leq l \leq d \\ d+1 \leq k \leq n}} (1-p^{-1}z_k/z_l)^{-1}.$$

\vspace{.1in}

\paragraph{\textbf{4.4. PBW theorem}} We finish this paragraph with the following important structural result due to Davison and Meinhardt~:

\begin{theorem}[Davison-Meinhardt, \cite{DM}]\label{T:DM} There exists a graded algebra filtration $\mathbb{Q} =F_0 \subseteq F_1 \subseteq \cdots$ of $H_*(T^*\mathcal{M}_Q)$ and an algebra isomorphism
\begin{equation}\label{E:PBW}
gr_{F_{\bullet}} (H_*(T^*\mathcal{M}_Q)) \simeq Sym(\bar{\mathfrak{n}}_Q[u])
\end{equation} 
where $\bar{\mathfrak{n}}_Q=\bigoplus_{\d \in \N^I} \bar{\mathfrak{n}}_{Q,\d}$ is a $\mathbb{N}^I \times \mathbb{N}$-graded vector space and $deg(u)=-2$. The same holds for $H_*(\Lambda_Q)$.
\end{theorem}

The filtration $F_{\bullet}$ is, essentially, the perverse filtration associated to the projection from the stack $T^*\mathcal{M}_Q$ to its coarse moduli space. As a direct corollary, $F_1 \simeq \bar{\mathfrak{n}}_Q[u]$ is equipped with the structure of an $\mathbb{N}$-graded Lie algebra; it is easy to see that it is the polynomial current algebra of an $\mathbb{N}$-graded Lie algebra $\bar{\mathfrak{n}}_Q$ (coined the \textit{BPS Lie algebra} in \cite{DM}). Loosely speaking, Theorem~\ref{T:DM} says that $H_*(T^*\mathcal{M}_Q)$ is (a filtered deformation of) the enveloping algebra $U(\bar{\mathfrak{n}}_Q[u])$, i.e. some kind of Yangian of $\bar{\mathfrak{n}}_Q$. Comparing graded dimension and using Theorem~\ref{T:BSV2} we deduce that for any $\d \in \N^I$
\begin{equation}\label{E:gradedimDM}
\text{dim}_{\Z}\;\bar{\mathfrak{n}}_{Q,\d}=A_{Q,\d}(t).
\end{equation}
 Moreover, when $Q$ has no edge loops, it can be shown that the degree zero Lie subalgebra $\bar{\mathfrak{n}}_Q[0]$ is isomorphic to the positive nilpotent subalgebra $\mathfrak{n}_Q$ of the Kac-Moody algebra $\mathfrak{g}_Q$. This implies at once both the positivity and the constant term conjectures for Kac polynomials, see Section~\textbf{2.2.}
The same reasoning also yields a proof of the nilpotent versions of the Kac conjectures (for an arbitrary quiver) of Section~\textbf{2.5.}, see \cite{Davison}.

\vspace{.1in}

At this point, the following conjecture appears inevitable~:

\begin{conjecture}\label{Conj:Hall=CoHall} The graded Lie algebra $\bar{\mathfrak{n}}_Q$ is isomorphic to the positive subalgebra $\widetilde{\mathfrak{n}}_Q$ of $\widetilde{\mathfrak{g}}_Q$.
\end{conjecture}
Notice that the definition of $\widetilde{\g}_Q$ involves the \textit{usual} Hall algebra of $Q$ (over all the finite fields $\mathbb{F}_q$), while that of $\bar{\mathfrak{n}}_Q$ involves the (two-dimensional) \textit{cohomological} Hall algebra of $Q$ and the complex geometry of $\text{Rep}_{\mathbb{C}}\, \Pi_Q$. Conjecture~\ref{Conj:Hall=CoHall}  would follow from the fact that $\bar{\mathfrak{n}}_Q$ is the positive half of some graded Borcherds algebra $\bar{\mathfrak{g}}_Q$.

\vspace{.1in}

\paragraph{\textbf{4.5.\;Action on Nakajima quiver varieties}} One important feature of the cohomological Hall algebras $H_*^T(T^*\mathcal{M}_Q), H^T_*(\Lambda_Q)$ is that they act, via some natural correspondences, on the cohomology of Nakajima quiver varieties.  Recall that the Nakajima quiver variety $\mathfrak{M}_Q(\v,\w)$ associated to a pair of dimension vectors $\v,\w \in \N^I$ is a smooth quasi-projective symplectic variety, which comes with a proper morphism $\pi~:\mathfrak{M}_Q(\v,\w) \to \mathfrak{M}_{Q,0}(\v,\w)$ to a (usually singular) affine variety. The morphism $\pi$ is an example of a symplectic resolution of singularities, of which quiver varieties provide one of the main sources. The quiver variety also comes with a canonical (in general singular) lagrangian subvariety $\mathfrak{L}(\v,\w)$ which, when the quiver has no edge loops, is the central fiber of $\pi$.
Examples include the Hilbert schemes of points on $\mathbb{A}^2$ or on Kleinian surfaces, the moduli spaces of instantons on theses same spaces, resolutions of Slodowy slices in nilpotent cones and many others (see e.g. \cite{SNaksurvey} for a survey of the theory of Nakajima quiver varieties). The following theorem was proved by Varagnolo in \cite{Varagnolo}, based on earlier work by Nakajima in the context of equivariant $K$-theory (see \cite{NakDuke98}, \cite{NakJAMS}).

\vspace{.1in}

\begin{theorem}[Varagnolo, \cite{Varagnolo}] Let $Q$ be without edge loops and $\w \in \N^I$. There is a geometric action of the Yangian $Y_h(\mathfrak{g}_Q)$ on $F_{\w}:=\bigoplus_{\v}H^{T\times G_{\w}}_*(\mathfrak{M}_Q(\v,\w))$, preserving the subspace $L_{\w}:=\bigoplus_{\v} H^{T\times G_{\w}}_*(\mathfrak{L}_Q(\v,\w))$. For $Q$ of finite type, $L_{\w}$ is isomorphic to the universal standard $Y_h(\mathfrak{g}_Q)$-module of highest weight $\w$.
\end{theorem}

In the above, the Yangian $Y_h(\mathfrak{g}_Q)$ is defined using Drinfeld's new realization, applied to an arbitrary Kac-Moody root system, see \cite{Varagnolo}. Its precise algebraic structure is only known for $Q$ of finite or affine type.

\vspace{.1in}

\begin{theorem}[S.-Vasserot, \cite{SVCOHA1}]\label{T:SVCOHA1} For any $Q$ and any $\w \in \N^I$ there is a geometric action of the cohomological Hall algebra $H_*^T(T^*\mathcal{M}_Q)$ on $F_\w$. The diagonal action of $H_*^T(T^*\mathcal{M}_Q)$ on $\prod_{\w} F_\w$ is faithful. The subalgebra $H_*^T(\Lambda_Q)$ preserves $L_{\w}$, which is a cyclic module.
\end{theorem}

\vspace{.1in}

For quivers with no edge loops, we recover Varagnolo's Yangian action by Theorem~\ref{T:generation}; It follows that there exists a surjective map $Y^+_h(\mathfrak{g}_Q)\otimes \mathbb{Q}(q_1,q_2) \to H^T_*(T^*\mathcal{M}_Q) \otimes \mathbb{Q}(q_1,q_2)$. This map is an isomorphism for $Q$ of finite or affine type. The action of $H^T_*(T^*\mathcal{M}_Q)$
is constructed by means of general Hecke correspondences; in fact one can view $H^T_*(T^*\mathcal{M}_Q)$ as the \textit{largest} algebra acting on $F_\w$ via Hecke correspondences.
Considering dual Hecke correspondences (or adjoint operators), one defines an opposite action of $H^T_*(T^*\mathcal{M}_Q)$; it is natural to expect that these two actions extend
to an action of some 'Drinfeld double' of $H^T_*(T^*\mathcal{M}_Q)$, but this remains to be worked out. It is also natural to expect that $L_{\w}$ is a universal or standard module
for $H^T_*(T^*\mathcal{M}_Q)$, as is suggested by Hausel's formula for the Poincar\'e polynomial of $\mathfrak{M}_Q(\v,\w)$ or $\mathfrak{L}_Q(\v,\w)$, which involves the (full) Kac polynomials, see \cite{HauselKac}, \cite{BSV}.  Theorem~\ref{T:SVCOHA1} has an obvious analog (with the same proof) for an arbitrary OBM theory (see \cite{YZ} for the construction
of the action).

\vspace{.1in}

\paragraph{\textbf{4.6.\;Relation to Maulik-Okounkov Yangians}} We finish this section by very briefly mentioning yet another (conjectural) construction of the Lie algebra $\widetilde{\g}_Q$, this time directly by means of the symplectic geometry of Nakajima quiver varieties. Using the theory of \textit{stable enveloppes} for $\mathbb{C}^*$-actions
on smooth symplectic varietes, Maulik and Okounkov constructed for any pair of dimension vectors $\w_1, \w_2$ a quantum $R$-matrix 
$$R_{\w,\w_2}(t)=1 + \frac{\hbar}{v^{-1}}  \mathbf{r}_{\w_1,\w_2} + O(v^{-2}) \in 
\text{End}(F_{\w_1} \otimes F_{\w_2})[[v^{-1}]].$$
 Applying the RTT formalism, we obtain a graded algebra $\mathbb{Y}_{Q}$ acting on all the spaces $F_\w$; similarly, from the classical $R$-matrices $\mathbf{r}_{\w_1,\w_2}$ we obtain a  $\mathbb{Z}$-graded Lie algebra $\mathbf{g}_Q$. Moreover, $\mathbb{Y}_Q$ is (up to some central elements) a filtered deformation of $U(\mathbf{g}_Q[u])$ (hence the name \textit{Maulik-Okounkov Yangian} see \cite{MO}). It can be shown that $\mathbf{g}_Q$ is a graded Borcherds algebra. The following conjecture was voiced by Okounkov (\cite{OkounkovTalk})~:
 
 \begin{conjecture}[Okounkov] For any $Q$ and any $\mathbf{d} \in \N^I$ we have $\text{dim}_{\mathbb{Z}} \;\mathbf{g}_Q[\d]=A_{Q,\d}(t)$.
 \end{conjecture}

The conjecture is known when $Q$ is of finite type. Comparing with (\ref{E:gradedimBKM}) and (\ref{E:gradedimDM}) leads to the following

\begin{conjecture} For any $Q$ we have $\mathbf{g}_Q\simeq \bar{\mathfrak{g}}_Q \simeq \widetilde{\mathfrak{g}}_Q$.\end{conjecture}

As a first step towards the above conjecture, we have

\begin{theorem}[S.-Vasserot, \cite{SVCOHA2}] For any $Q$ there is a canonical embedding $H^T_*(\Lambda_Q) \to \mathbb{Y}^+_Q$, compatible with the respective actions on $\prod_\w F_\w$
of $H^T_*(\Lambda_Q)$ and $\mathbb{Y}_Q$.
\end{theorem}

\vspace{.1in}

All together, we see that there are conjecturally (at least) \textit{three} different incarnations of the \textit{same} graded Borcherds Lie algebra : as a 'generic form' of the full Hall algebra of the category of representations of $Q$ over finite fields, as a cohomological Hall algebra of the complex (singular) stack $T^*\mathcal{M}_Q$, or as an algebra
acting on the cohomology of Nakajima quiver varieties via R-matrices constructed by means of symplectic geometry. The precise structure of this graded Borcherds algebra, in particular
its Cartan matrix (determined by the polynomials $C_{Q,\d}^{\text{abs}}(t)$ counting the dimensions of the spaces of 'absolutely' cuspidal functions for $Q$) remains however very mysterious. From the symplectic geometry perspective, one might expect the polynomials $C_Q^{\text{abs}}(t)$ to be related to some motive inside the Nakajima quiver variety, but there is no conjectural construction of such a motive that we know of.

\bigskip

In the remainder of this paper, we shift gears and consider categories of coherent sheaves on smooth projective curves instead of representations of quivers; motivated by the analogy with quivers, we will describe the spherical Hall algebra, full Hall algebra and $2d$-cohomological Hall algebras (!) of a curve, as well as a geometric interpretation of Kac polynomials
in terms of Higgs bundles. We finish with some speculation about the Lie theoretic structures that we believe are lurking in the background.

 \vspace{.1in}
 
\section{Hall algebras of curves and shuffle algebras}
%3 pages
 \medskip
 
 \paragraph{\textbf{5.1. Notations }} We fix an integer $g \geq 0$ and a smooth, geometrically connected, projective curve $X$ of genus $g$ over a field $\k$. We denote by $\text{Coh}(X)$ the category of coherent sheaves on $X$ and by $\mathcal{M}_{X,r,d}$ the stack of coherent sheaves of rank $r$ and degree $d$. It is a smooth stack, locally of finite type. The Euler form is given by the Riemann-Roch formula
 $$\langle \mathcal{F},\mathcal{G}\rangle= (1-g) rk(\mathcal{F})rk(\mathcal{G}) + (rk(\mathcal{F})deg(\mathcal{G})-rk(\mathcal{G})deg(\mathcal{F})).$$
 Set $\mathcal{M}_X=\bigsqcup_{r,d} \mathcal{M}_{X,r,d}$. 
 
 \vspace{.1in}
 
 \paragraph{\textbf{5.2. Ringel-Hall algebra of a curve}} Let us now assume that $\k$ is a finite field, and define
$$\H_X=\bigoplus_{r,d} \H_X[r,d], \qquad \H_X[r,d]:=\text{Fun}(\mathcal{M}_{X,r,d}(\k) ,\C).$$
In an entirely similar fashion to (\ref{E:conv}), there is a convolution diagram
\begin{equation}\label{E:conv2}
\xymatrix{ \mathcal{M}_X \times \mathcal{M}_X & \widetilde{\mathcal{M}}_X \ar[l]_-q \ar[r]^-p & \mathcal{M}_X}
\end{equation}
where $\widetilde{\mathcal{M}}_X$ is the stack parametrizing short exact sequences 
$$\xymatrix{0 \ar[r] & \mathcal{F} \ar[r] & \mathcal{H} \ar[r] & \mathcal{G} \ar[r]& 0}$$
in $\text{Coh}(X)$; The map $p$ is proper, while the map $q$ is again a stack vector bundle, whose restriction to 
 $\mathcal{M}_{X,r,d} \times \mathcal{M}_{X,r',d'}$ is of rank $-\langle (r',d'),(r,d)\rangle$. 
Setting $v=(\#\k)^{\frac{1}{2}}$, $\mathbf{K}=\C[k_{(0,1)}^{\pm 1}, k_{(1,0)}^{\pm 1}]$ we equip, using (\ref{E:prodhall}) and (\ref{E:prodhall2}), $\widetilde{\H}_X:=\H_X \otimes \mathbf{K}$ with the structure of a $\Z^2$-graded bialgebra called the \textit{Hall algebra} of $X$. It carries a nondegenerate Hopf pairing by (\ref{E:prodhall3}). Finally,
$\H_X \simeq \H_X^{\text{bun}} \ltimes \H_X^0$ where $\H_X^{\text{bun}}, \H_X^{0}$ are the Hall subalgebras of vector bundles, resp. torsion sheaves.
The Hall algebra of a curve was first considered by Kapranov in \cite{KapEis}, who established the (direct)
dictionary between $\H_X^{\text{bun}}$ and the space of automorphic forms\footnote{one nice feature of the Hall algebra is that it incorporates the algebra of Hecke operators, as the sub Hopf algebra of torsion sheaves $\H_X^0$} for the groups $GL(n, \mathbb{A}_X)$, together with the operations of parabolic induction (Eisenstein series) and restriction (constant term map), and who observed a striking analogy between $\widetilde{\H}_X$ and quantum loop groups.

Let $\mathbf{H}_X^{\text{sph}}$ (resp. $\H_X^{\text{sph,bun}}$) be the subalgebra of $\mathbf{H}_X$ generated by the constant functions on $\mathcal{M}_{X,r,d}$ for $r=0,1$ (resp. $r=1$) and $d \in \Z$. From the point of view of automorphic forms, $\H_X^{\text{sph,bun}}$ is the space spanned by all components of Eisenstein series induced from trivial (i.e constant) automorphic forms for the torus (together with a suitable space of Hecke operators in the case of $\H_X^{\text{sph}}$). One can show that $\widetilde{\mathbf{H}}^{\text{sph}}_X:={\mathbf{H}}^{\text{sph}}_X \otimes \mathbf{K}$ is a self-dual sub-Hopf algebra of $\widetilde{\H}_X$ which contains the characteristic functions of any Harder-Narasimhan strata (see \cite{SKorea}). Moreover, contrary to $\widetilde{\H}_X$ which depends strongly on the fine arithmetic structure of $X$, $\widetilde{\H}_X^{\text{sph}}$ only depends on the Weil numbers of $X$ and admits an $R_g$-\textit{rational form}, i.e. there exists
a torsion-free $R_g$-Hopf algebra $\widetilde{\H}_{\Sigma_g}^{\text{sph}}$ such that for any smooth projective curve $X$ defined over a finite field, $\widetilde{\H}_{\Sigma_g}^{\text{sph}} \otimes_{R_g} \C_X \simeq \widetilde{\H}_X^{\text{sph}}$, where 
$\C_X$ is the $R_g$-module corresponding to the evaluation morphism $R_g \to \overline{\mathbb{\Q}_l} \simeq \C$, $f \mapsto f(Fr_x)$. We view $\widetilde{\H}_{\Sigma_g}^{\text{sph}}$ as some kind of (half) quantum group which depends on $dim(T_g)=g+1$ quantum parameters, associated to curves of genus $g$. The full quantum group is, as before, obtained by the Drinfeld double procedure. Hall algebras in genus $0$ and $1$ are already very interesting~:

\begin{theorem}[Kapranov, \cite{KapEis}, Baumann-Kassel, \cite{BK}] The Drinfeld double $\mathbf{D}\widetilde{\H}^{\text{sph}}_{\Sigma_0}$ is isomorphic to the quantum affine algebra $U_v(\widehat{\mathfrak{sl}_2})$.
\end{theorem}

\begin{theorem}[Burban-S., \cite{BS}, S.-Vasserot, \cite{SVCompos}] The Drinfeld double $\mathbf{D}\widetilde{\H}^{\text{sph}}_{\Sigma_1}$ is isomorphic to the spherical double affine Hecke algebra
$\mathbf{S}\ddot{\mathbf{H}}_{q,t}(GL_\infty)$ of type $GL(\infty)$.
\end{theorem}
 
 In the above two cases, the structure of the Hall algebra $\widetilde{\H}^{\text{sph}}_{\Sigma_g}$ is rather well understood : it has a PBW-type basis as well as a canonical basis constructed from simple perverse sheaves (Eisenstein sheaves) on the stacks $\mathcal{M}_X$. The spherical Hall algebra $\widetilde{\H}^{\text{sph}}_{\Sigma_1}$ --also called the \textit{elliptic Hall algebra}-- has found a surprising number of applications in representation theory of Cherednik algebras, low-dimensional topology and knot theory (e.g. \cite{MS} , \cite{GorkskyNegut}), algebraic geometry and mathematical physics of the instanton spaces on $\mathbb{A}^2$ (e.g. \cite{SVDuke}, \cite{SVIHES}, \cite{NegutAGT}), combinatorics of Macdonald polynomials (e.g. \cite{Bergeron}, \cite{DiFrancesco}), categorification (e.g. \cite{CLLSS}), etc. The elliptic Hall algebra (or close variants thereof) has independently appeared in the work of Miki, \cite{Mikii}, Ding-Iohara, \cite{Ding} and Feigin and his collaborators (see e.g., \cite{Feigin}). 

\vspace{.1in}
 
\paragraph{\textbf{5.3.\;Shuffle algebra presentation}} Although the structure of $\widetilde{\H}^{\text{sph}}_{\Sigma_g}$ for $g >1$ is much less well understood, we always have a purely algebraic model of $\widetilde{\H}_{\Sigma_g}^{\text{sph}}$ at our disposal, once again\footnote{There is nothing surprising in the ubiquity of shuffle algebras~: any finitely generated, $\mathbb{N}$-graded self-dual Hopf algebra has a realization as a shuffle algebra} in the guise of a shuffle algebra. More precisely, let
$$\zeta_{\Sigma_g}(z):=\frac{\prod_{i=1}^g (1-\eta_iz)}{(1-z)(1-qz)} \in R_g(z)$$
be the 'generic' zeta function of a curve of genus $g$ and put $\zeta'_{\Sigma_g}(z)=(1-qz)(1-qz^{-1})\zeta_{\Sigma_g}(z)$. Consider the shuffle algebra
$$Sh_{\Sigma_g}:=\bigoplus_d R_g[z_1^{\pm 1}, \ldots, z_d^{\pm 1}]^{\mathfrak{S}_d}$$
with multiplication
\begin{equation}\label{E:shuffle2}
(f *g)(\underline{z}_{[1,n]})=\sum_{ \sigma \in Sh_{d,e}} \sigma \left[ K_{d,e}(\underline{z}_{[1,n]}) f(\underline{z}_{[1,d]})g(\underline{z}_{[d+1,n]})\right]
\end{equation}
where $K_{d,e}(\underline{z}_{[1,n]})=\prod_{1\leq l\leq d<k\leq n} \zeta'_{\Sigma_g}(z_l/z_k)$, for any $d,e$ and $n=d+e$.

\medskip

\begin{theorem}[S.-Vasserot, \cite{SVMathann}]\label{T:geomlang} The assignment $1_{Pic^d} \mapsto z^d \in Sh_{\Sigma_g}[1]$ extends to an $R_g$-algebra embedding $\Psi~:\H_{\Sigma_g}^{\text{sph,bun}} \to
Sh_{\Sigma_g}$.
\end{theorem}

The map $\Psi$ is essentially the iterated coproduct $\Delta^{(r)}~: \H^{\text{sph}}_{\Sigma_g}[r] \to (\H^{\text{sph}}_{\Sigma_g}[1] )^{\otimes r}$. In the language of automorphic forms,
Theorem~\ref{T:geomlang} amounts to the Langlands formula for the constant term of Eisenstein series, or equivalently to the Gindikin-Karpelevich formula. One should not be deceived by the apparent simplicity of the shuffle algebra description for $\H^{\text{sph,bun}}_X$. In particular, $\H^{\text{sph,bun}}_{\Sigma_g}$ is \textit{not} free over $R_g$ and the relations satisfied by the generators $1_{Pic^d}$ of $\H_X^{\text{sph,bun}}$ \textit{do} depend on the arithmetic of the Weil numbers $\sigma_1, \ldots, \sigma_{2g}$ of $X$ (more precisely, on the $\mathbb{Z}$-linear dependence relation between $log(\sigma_1), \ldots log(\sigma_{2g})$, the so-called \textit{wheel relations}).

\vspace{.1in}

\paragraph{\textbf{5.4. Geometric Langlands isomorphism}} Comparing (\ref{E:shuffle1}) and (\ref{E:shuffle2}) we immediately see that, up to the identification $p=q^{-1}, \xi_i= \eta_i^{-1}$ we have $Sh_{\Sigma_g} \simeq Sh_{S_g}^{K^{T_g}}$. This implies that there is an isomorphism
\begin{equation}\label{E:geomlang1}
\Phi~: \H_{\Sigma_g}^{\text{sph,bun}} \stackrel{\sim}{\longrightarrow} \mathbf{K}_{S_g}^{\text{sph}, T_g}
\end{equation}
where $\mathbf{K}_{S_g}^{\text{sph}, T_g}$ is the subalgebra of $K^{T_g}(T^*\mathcal{M}_{S_g})$ generated by its rank one component $K^{T_g}(T^*\mathcal{M}_{S_g,1})$. The
isomorphism (\ref{E:geomlang1}) between the (spherical) Hall algebra of $X$ and the (spherical) K-theoretical Hall algebra of the $S_g$ quiver should be viewed as an incarnation, at the level of Grothendieck groups, of a (linearized) form of geometric Langlands correspondence. Indeed, the Langlands philosophy predicts an equivalence
$$Coh(LocSys_r(X_{\C})) \simeq D\hspace{-.04in}-\hspace{-.04in}mod(Bun_r(X_\C))$$
between suitable (infinity) categories of coherent sheaves on the moduli stack of $GL_r$-local systems on $X_\C$ and $D$-modules on the stack of $GL_r$-bundles on $X_\C$. On one hand, the formal neighborhood in $LocSys_{r}(X_\C)$ of the trivial local system may be linearized as the formal neighborhood of $0$ in 
$\{(x_1, \ldots, x_g, y_1, \ldots, y_g) \in \mathfrak{gl}_r(\C)^{2g}\;|\; \sum_i [x_i,y_i]=0\}/ GL_r(\C) \simeq T^*\mathcal{M}_{S_g,r}$; on the other hand, $\H^{sph,bun}_{\Sigma_g}$ may be lifted to a category of holonomic D-modules (perverse sheaves) on $Bun_r(X_{\C})$. We refer to \cite{SVMathann} for a more detailed discussion.

\vspace{.1in}

\paragraph{\textbf{5.5.\;Variations}} There are several interesting variants of the above constructions and results : one can consider categories of $D$-parabolic coherent sheaves
(see \cite{SDukeI}, \cite{Lin}); this yields, for instance, some quantum affine or toroidal algebras. One can also consider arithmetic analogs the Hall algebra, replacing the abelian category of coherent sheaves on a curve $X$ over a finite field by coherent sheaves or vector bundles (in the sense of Arakelov geometry) over $\overline{Spec(\mathcal{O}_K)}$, where $K$ is a number field. The case of $K=\mathbb{Q}$ is discussed in \cite{KSV}, where the spherical Hall algebra is described as some analytic shuffle algebra with kernel
given by the Riemann zeta function $\zeta(z)$.

\medskip
 
%\section{The elliptic Hall algebra}
%%2 pages
%\medskip
%
%\paragraph{\textbf{6.1.}\textit{\;Moduli stacks of coherent sheaves on elliptic curves.}} just some notations and Atiyah's theorem
%
%\paragraph{\textbf{6.2.}\textit{\;A diagrammatic presentation.}} including invariance under $SL(2,\Z)$ and derived auto-equivalences. Shuffle algebra presentation; Mention Dragos theorem. The Lie algebra here is toroidal $gl(1)$
%
%\paragraph{\textbf{6.3.}\textit{\;Spherical Cherednik algebra and Macdonald polynomials}} including construction as Eisenstein sheaves restricted to semistable locus; Fourier = Fourier-Mukai transform; symmetric invariants; pure combinatorics -- mention Bergeron business.
%
%\paragraph{\textbf{6.4.}\textit{\;Dual picture : action on the K-theory of instanton spaces}} generalize Nakajima's picture. Fock space representation for Hilbert scheme case. Mention Tsymbaliuk and Feigin's works. Mention AGT ! here
%
%\paragraph{\textbf{6.5.}\textit{\;Other incarnations~: Torus knots, skein algebras, $q$-Heisenberg category}} Samuelson and 
%link to torus knot invariants ----need more here !!!
%
% \medskip
 
\vspace{.1in}
 
\section{Counting Cuspidals and Cohomological Hall algebras of curves}
%1,5pages
\medskip

\paragraph{\textbf{6.1.\;Counting cuspidals and the full Hall algebra of a curve}} In Sections~\textbf{5.2., 5.3.} we considered the spherical Hall algebra $\H_X^{\text{sph}}$ for which we provided a shuffle algebra description involving the zeta function $\zeta_X(z)$. What about the whole Hall algebra $\H_X$ ? As shown in \cite{KSV2}, $\H^{bun}_X$ admits a shuffle description as well, but it is much less explicit than for $\H_X^{\text{sph}}$. Recall that an element $f \in \H_X^{\text{bun}}[r,d]$ is \textit{cuspidal} if it is quasi-primitive, i.e. if $$\Delta(f) \in f \otimes 1 +k_{r,d} \otimes f +\H_X^0 \otimes \widetilde{\H}_X[r].$$ The algebra $\H^{\text{bun}}_X$ is generated by the spaces of cuspidal elements
$\H_X^{\text{cusp}}=\bigoplus_{r,d} \H_X^{\text{cusp}}[r,d]$, and $\dim\;\H_X^{\text{cusp}}[r,d] < \infty$ for all $r,d$.  The function field Langlands program (\cite{Lafforgue}) sets up a
correspondence $\chi \mapsto f_{\chi}$ between characters $\chi~: \H_X^0 \to \C$ associated to rank $r$ irreducible local systems on $X$ and cuspidal Hecke eigenfunctions $f_\chi = \sum_d f_{\chi,d} \in \prod_d \H_{X}^{\text{cusp}}[r,d]$. The shuffle algebra description of $\H_X$ is as follows~: we have a family of variables $z_{\chi,i}$, $i \in \N$ for \textit{each} cuspidal Hecke eigenform $f_\chi$ (up to $\mathbb{G}_m$-twist) and the shuffle kernels involve the Rankin-Selberg $L$-functions $L(\chi,\chi',z)$ of pairs of characters $\chi, \chi' $ in place of $\zeta_X(z)$ (see \cite{Dragos} for a full treatment when $g=1$). In principle, one could try, using the PBW theorem and arguing as in the case of quivers (see Section~\textbf{3.2.}), to deduce from the above shuffle description of $\H_X$ an expression for the dimensions of the spaces of cuspidal functions $\H_X^{\text{cusp}}[r,d]$ (or better, absolutely cuspidal functions $\H_X^{\text{abs. cusp}}[r,d]$) in terms of the Kac polynomials $A_{g,r}$ . Very recently H. Yu managed --by other means\footnote{namely, using the Arthur-Selberg trace formula}-- to compute the dimension of $\H_X^{\text{abs. cusp}}[r,d]$ directly~:

\begin{theorem}[Yu, \cite{Yu}]\label{T:Yu} For any $g,r$ there exists a (unique) polynomial $C^{\text{abs}}_{g,r} \in R_g$ such that for any smooth projective curve $X$ of genus $g$ defined over
a finite field, $\dim\; \H_X^{\text{abs. cusp}}[r,d]= C^{abs}_{g,r}(Fr_X)$.
\end{theorem}

This generalizes a famous result of Drinfeld (for $r=2$, \cite{Drinfeld}) and proves a conjecture of Deligne (\cite{Deligne}) and Kontsevich (\cite{Kontsevich}). Interestingly, the polynomial
$C^{\text{abs}}_{g,r}$ is explicit : Yu expresses it in terms of the numbers of rational points of the moduli spaces of stable Higgs bundles over finite fields and hence, by Theorem~\ref{T:last} below,
in terms of the Kac polynomials $A_{g,r}$ ! This strengthens our belief that the structure of $\H_X$ as an associative algebra is nice enough that it should have a character formula similar to that of Borcherds algebras.
For instance, we have
$$A_{g,1}(F)=C^{\text{abs}}_{g,1}(F), \qquad A_{g,2}(F)=C^{\text{abs}}_{g,2}(F) + (g-1)C_{g,1}^{\text{abs}}(F)^2+C_{g,1}^{\text{abs}}(F)$$
\begin{equation*}
\begin{split}
A_{g,3}(F)=&C^{\text{abs}}_{g,3}(F) + (g-1)C_{g,1}^{\text{abs}}(F) \left\{4C^{\text{abs}}_{g,2}(F) + C^{\text{abs}}_{g,1}(F^2) + 2(g-1)C^{\text{abs}}_{g,1}(F)^2\right\}\\
&+ 4(g-1) C^{\text{abs}}_{g,1}(F)^2+C^{\text{abs}}_{g,1}(F)
\end{split}
\end{equation*}
where $F=Fr_X$.

\begin{conjecture}\label{C:kacpolcat2} For any $g, r$ there exists a (non-virtual) $GSp(2g,\overline{\mathbb{Q}_l})$-representation $\mathbb{C}^{\text{abs}}_{g,r}$ such that $C^{\text{abs}}_{g,r}=\tau ( ch(\mathbb{C}^{\text{abs}}_{g,r}))$.
\end{conjecture}

\vspace{.1in}

\paragraph{\textbf{6.2.\;Cohomological Hall algebra of Higgs sheaves.}} We take $\k=\C$ here. The (underived) cotangent stack $T^*\mathcal{M}_X=\bigsqcup_{r,d} T^*\mathcal{M}_{X,r,d}$ is identified  with the stack of Higgs sheaves $\mathcal{H}iggs_X=\bigsqcup_{r,d} \mathcal{H}iggs_{X,r,d}$, which parametrizes pairs $(\mathcal{F}, \theta)$ with $\mathcal{F} \in \text{Coh}(X)$ and $\theta \in Hom_{\mathcal{O}_X}(\mathcal{F}, \mathcal{F} \otimes \Omega_X)$. The \textit{global nilpotent cone} $\Lambda_{X} = \bigsqcup_{r,d} \Lambda_{X,r,d}$ is the closed Lagrangian substack whose objects are Higgs sheaves $(\mathcal{F},\theta)$ for which $\theta$ is nilpotent. Both stacks are singular, locally of finite type and have infinitely many irreducible components. The stack $\Lambda_{X,r,d}$ is slightly better behaved since all irreducible components are of dimension $(g-1)r^2$. We refer to \cite{TristanNilp} for an explicit description of these irreducible components. The torus $T=\mathbb{G}_m$
acts on $\mathcal{H}iggs_X$ and $\Lambda_X$ by $t \cdot (\mathcal{F}, \theta)=(\mathcal{F}, t\theta)$.

\begin{theorem}[Sala-S., \cite{SalaS}, Minets, \cite{Minets}] The Borel-Moore homology spaces $H^T_*(\mathcal{H}iggs_X)$ and $H^T_*(\Lambda_X)$ carry natural $\mathbb{Z} \times \mathbb{Z}^2$-graded associative algebra structures. Moreover, the direct image morphism $i_* :H^T_*(\Lambda_X) \to H_*^T(\mathcal{H}iggs_X)$ is an algebra homomorphism.
\end{theorem}

The definition of the algebra structure roughly follows the same strategy as for quivers~: working in local charts, we use the construction of $\mathcal{H}iggs_X$ as
a symplectic reduction to embed everything into some smooth moduli stacks. As before, the morphism $i_*$ becomes invertible after localizing with respect to $H^*_T(pt)$.
There is an embedding of $\mathcal{M}_X$ in $\Lambda_X$ as the zero section of the projection $p: \mathcal{H}iggs_X \to \mathcal{M}_X$; its image is an irreducible component of $\Lambda_X$.

\begin{theorem}[Sala-S., \cite{SalaS}] The algebra $H_*^T(\Lambda_X)$ is generated by the collection of subspaces $H_*^T(\mathcal{M}_{X,r,d})$ for $(r,d) \in \Z^2$.
\end{theorem}

What about a shuffle algebra description of $H_*^T(\Lambda_X)$ ? The cohomology ring $\mathbb{H}_{r,d}=H^*(\mathcal{M}_{X,r,d})$ acts on
$H^T_*(\Lambda_{X,r,d})$ by $c \cdot h=p^*(c) \cap h$.
By Heinloth's extension of the Atiyah-Bott theorem \cite{HeinlothAB}, $\mathbb{H}_{0,d}=S^d(H^*(X)[z])$ for $d \geq 1$ while for $r \geq 1, d \in \Z$, $\mathbb{H}_{r,d}=\Q[c_{i,\gamma}(\mathcal{E}_{r,d})\;|\; i \geq 1, \gamma]$ is a polynomial algebra in the K\"unneth components $c_{i,\gamma}(\mathcal{E}_{r,d})$ of the Chern classes
of the tautological sheaf $\mathcal{E}_{r,d}$ over $\mathcal{M}_{X,r,d} \times X$; here $\gamma$ runs over a basis $\{1, a_1, b_1, \ldots, a_g, b_g, \varpi\}$ of $H^*(X)$.

\begin{theorem}[Sala-S., \cite{SalaS}, Minets, \cite{Minets}] The $\mathbb{H}_{r,d}$-module $H^T_*(\Lambda_{X,r,d})$ is torsion-free and of generic rank one.
\end{theorem}

This opens up the possibility to construct a (shuffle) algebra structure on $\bigoplus_{r,d} \mathbb{H}_{r,d}$, but this has so far only be achieved for the subalgebra
of torsion sheaves, see \cite{Minets}. Very similar shuffle algebras occur, not surprisingly, as operators on the cohomology of moduli spaces of semistable sheaves on smooth surfaces, see \cite{Negut}.

\vspace{.1in}

\section{Kac polynomials and Poincar\'e polynomials of moduli of stable Higgs bundles}
%1 page

 \medskip
 
\paragraph{\textbf{7.1. Moduli spaces of stable Higgs bundles}} Recall that a Higgs sheaf $(V,\theta)$ on $X$ is \textit{semistable} if for any subsheaf $W \subset V$ such that
$\theta(W) \subseteq W \otimes \Omega_X$ we have $\mu(W) \leq \mu(V)$, where $\mu(\mathcal{F})=deg(\mathcal{F})/rk(\mathcal{F})$ is the usual slope function. 
Replacing $\leq$ by $<$ we obtain the definition of a \textit{stable} Higgs sheaf. The open substack $\mathcal{H}iggs^{st}_{X,r,d} \subset \mathcal{H}iggs_{X,r,d}$ of stable sheaves is a $\mathbb{G}_m$-gerbe over a smooth quasi-projective symplectic variety $Higgs^{st}_{r,d}$. On the contrary, the open substack $\mathcal{H}iggs^{ss}_{r,d} \subset \mathcal{H}iggs_{X,r,d}$ of semistable sheaves is singular as soon as $gcd(r,d) >1$; when $gcd(r,d)=1$, $\mathcal{H}iggs^{ss}_{X,r,d}=\mathcal{H}iggs^{st}_{X,r,d}$. The variety $Higgs_{X,r,d}^{st}$ has played a fundamental role in algebraic geometry, in the theory of integrable systems, in the geometric Langlands program, in the theory of automorphic forms, and is still the focus of intensive research; we refer to e.g. \cite{Hauselsurvey} for a survey of the many important conjectures in the subject.

 \vspace{.1in}
 
\paragraph{\textbf{7.2. Poincar\'e polynomials and Kac polynomials}} Assume that $\k=\mathbb{F}_q$. The following result provides a very vivid geometric interpretation of the Kac polynomial $A_{g,r}$.

\begin{theorem}[S., \cite{Sannals}, S.-Mozgovoy, \cite{MozSchif}]\label{T:last} For any $r,d$ with $gcd(r,d)=1$ we have 
$$\#Higgs_{X,r,d}^{st}(\mathbb{F}_q)=q^{(g-1)r^2+1}A_{g,r}(Fr_X).$$
\end{theorem}

Together with some purity argument and Mellit's simplification of the explicit formula for $A_{g,r}$, this proves  the conjecture of Hausel and Rodriguez-Villegas \cite{HRVInvent} for the Poincar\'e polynomial of $Higgs^{st}_{X,r,d}$ (for $\k=\C$ or $\mathbb{F}_q$). There are two proofs of Theorem~\ref{T:last} : one is based on a deformation argument to relate directly the point count of $Higgs_{X,r,d}(\mathbb{F}_q)$ to the count of indecomposable sheaves; the other uses Hall-theoretic techniques (for the category of positive Higgs sheaves). Both have been generalized or put in a broader context, see \cite{Ginzburg}, \cite{Fedorov}.

\medskip

Concerning the stable global nilpotent cone $\Lambda_{X,r,d}^{st}=\Lambda_{X,r,d} \cap Higgs_{X,r,d}^{st}$ we have the following interpretation of the constant term of the Kac polynomial ~:

\begin{theorem}[S., \cite{Sannals}] For any $(r,d)$ with $gcd(r,d)=1$ we have $\# Irr(\Lambda^{st}_{X,r,d})=A_{g,r,d}(0)$.
\end{theorem}

Recall that by Proposition~\ref{P:strangerelation} we have $A_{g,r}(0)=A_{S_g,r}(1)$; this suggests the existence of some natural partition of $Irr(\Lambda_{X,r,d}^{st})$, but the geometric meaning of such a partition is unclear to us.

\vspace{.1in}

\paragraph{\textbf{7.3. Donaldson-Thomas invariants and Kac polynomials}} What about \textit{non} coprime $(r,d)$~? In this case, it is still possible to
perform the (orbifold) point count of the stack $\mathcal{H}iggs_{X,r,d}(\k)$ when $\k=\mathbb{F}_q$; this point count is best expressed in terms of the Donaldson-Thomas invariants $\Omega_{X,r,d}$ which are defined by the following generating series~:
$$\forall \; \nu \in \mathbb{Q}, \qquad \sum_{\frac{d}{r}=\nu} \frac{\Omega_{X,r,d}}{q-1}w^rz^d := \text{Log}\left(\sum_{\frac{d}{r}=\nu} q^{(1-g)r^2} \#(\mathcal{H}iggs_{X,r,d}^{ss}(\mathbb{F}_q))w^rz^d\right).$$
 
 \begin{theorem}[S.-Mozgovoy, \cite{MozSchif}] For any $r,d$ we have $\Omega_{X,r,d}=qA_{g,r}(Fr_X)$.
\end{theorem}

If $gcd(r,d)=1$ then $\Omega_{X,r,d}=(q-1)q^{(1-g)r^2}\#(\mathcal{H}iggs_{X,r,d}^{st}(\mathbb{F}_q))=q^{(1-g)r^2}\#Higgs_{X,r,d}^{st}(\mathbb{F}_q)$, so that we recover Theorem~\ref{T:last}.

\vspace{.1in}

\section{Delirium Tremens :\;a hierarchy of Lie algebras}

\medskip

We conclude this survey with some wild speculations concerning potential Lie algebras associated to curves, rather than to quivers

\medskip

\paragraph{\textbf{8.1. Lie algebras from curves ?}} Following the analogy with quivers,
it is natural to expect the existence of a family of $\Z^2$-graded complex Lie algebras $\mathfrak{g}_g = \bigoplus_{r,d} \mathfrak{g}_{g}[r,d]$ such that $\H^{\text{sph}}_g$ is a $(g+1)$-quantum parameter deformation of $U^+(\mathfrak{g}_g)$, and 
$$\dim\;\mathfrak{g}_g[r,d]=A_{g,r}(0,\ldots, 0)$$ for any $r,d$. This Lie algebra $\g_g$
would be a curve analog of the Kac-Moody algebra $\g_Q$ associated to a quiver $Q$ (or its variant $\g^B_Q$ if $Q$ has edge loops). What about the the analog of the graded Borcherds algebra $\widetilde{\g}_Q$ ? Because the grading in the context of curves is by the character ring of $GSp(2g, \overline{\mathbb{Q}_l})$ rather than by $\Z$, it seems natural to expect the existence of a Lie algebra $\widetilde{\g}_g$ \textit{in the tensor category of finite-dimensional $GSp(2g,\overline{\mathbb{Q}_l})$-modules}, with $\g_g$ being identified with the sub-Lie algebra corresponding to the tensor subcategory of trivial representations (of arbitrary rank). Moreover, we should have
$$\widetilde{\g}_g[r,d] = \mathbb{A}_{g,r} \in GSp(2g,\overline{\mathbb{Q}_l})-\text{mod},$$ where $\mathbb{A}_{g,r}$ is as in Conjecture~\ref{C:kacpolcat}, and the cuspidal (or 'simple root vectors') of $\widetilde{\g}_g[r,d]$ should form a subrepresentation isomorphic to $\mathbb{C}^{\text{abs}}_{g,r}$, where $\mathbb{C}^{\text{abs}}_{g,r}$ is as in Conjecture~\ref{C:kacpolcat2}.

\medskip

 Although we do not have any clue at the moment as to what $\widetilde{\g}_g$ could be, the Langlands isomorphism (\ref{E:geomlang1}) provides us with a very good guess concerning $\g_g$. Namely, it is expected that the (spherical) K-theoretical Hall algebra $K^{T}(T^*\mathcal{M}_{Q})$-- by analogy with the (spherical) cohomological Hall algebra $H^{T}_*(T^*\mathcal{M}_{Q})$-- is a deformation of $U^+(\widetilde{\g}_{Q}[u^{\pm 1}])$. This strongly suggests that, at least as a vector space, $\g_g \simeq \widetilde{\g}_{S_g}[u^{\pm 1}]$. Note that the $\N$-grading of $\widetilde{\g}_{S_g}$ gets lost in the process since there is no obvious grading in the K-theoretical Hall algebra. Taking graded dimensions, we obtain the equality $A_{S_g,r}(1)=A_{g,r}(0,\ldots, 0)$ of Proposition~\ref{P:strangerelation}. Of course, this is not a proof but rather a conceptual explanation of this equality. We summarize this in the chain of inclusions of Lie algebras
$$\g_{S_g} \subseteq \tilde{\g}_{S_g} \subset \tilde{\g}_{S_g}[t,t^{-1}]\simeq \g_g \subset \tilde{\g}_g \subset  \tilde{\g}_g[t,t^{-1}].$$

\medskip

\noindent
\textit{Examples.} i) Suppose $g=0$. Then $\widetilde{\g}_{S_0} =\g_{S_0}= \mathfrak{sl}_2$ and $\g_0 \simeq \widehat{\mathfrak{sl}_2}$, while it is natural to expect that $\widetilde{\g}_0 = (\mathfrak{sl}_2 \oplus K)[u^{\pm 1}]$, where $K$ is a one dimensional central extension, placed in degree one. Note that we have $A_{0,0}=q+1$, $A_{0,1}=1$ and $A_{0,r}=0$ for $r>1$.\\
ii) Suppose $g=1$. Then $\widetilde{\g}_{S_1} \simeq \g_{S_1} = \overline{\mathbb{Q}_l}[s^{ \pm 1}]$ is the Heisenberg algebra, and $\g_1 = \overline{\mathbb{Q}_l}[s^{ \pm 1},t^{\pm 1}]\oplus K_1 \oplus K_2$, where $K_1\oplus K_2$ is a two-dimensional central extension. The Lie algebra structure is not the obvious one however,
but rather a central extension of the Lie bracket $[s^rt^d,s^nt^m]=(rm-dn)s^{r+n}t^{d+m}$ (see \cite{SVIHES}, App. F. for the case of the Yangian).

\vspace{.1in}

\paragraph{\textbf{8.2. Summary of Hall algebras, their corresponding Lie algebras and Kac polynomials.}} We conclude this survey with the following table, containing
our heuristics.

\vspace{.2in}

{\center{
\begin{tabular}{|c|c|c|}
  \hline
  Type of Hall Algebras & Quivers & Curves \\
  \hline
   $K$-theoretic Hall algebra $K^{sph,T}(T^*\mathcal{M})$& $\widetilde{\g}_Q[u^{\pm 1}]$; $A_Q(t)\delta(t)$  & $\widetilde{\g}_g[u]$; $A_g(\sigma_1, \ldots, \sigma_{2g})\delta(t)$ \\
  \hline
  Cohomological Hall algebra $H^T_*(T^*\mathcal{M})$&$\widetilde{\g}_Q[u]$; $A_Q(t)/(1-t)$  & $\widetilde{\g}_g[u]$; $A_g(\sigma_1, \ldots, \sigma_{2g})/(1-t)$ \\
  \hline
  Hall algebra $\H$ & $\widetilde{\g}_Q$; $A_Q(t)$ & $\widetilde{\g}_g$; $A_g(\sigma_1, \ldots, \sigma_{2g})$  \\
  \hline
  Spherical Hall algebra $\H^{\text{sph}}$ & $\g_Q$; $A_{Q}(0)$ & $\g_g$; $A_{g}(0,\ldots, 0)$ \\
  \hline

\end{tabular}}}

\vspace{.3in}

\paragraph{\textbf{Acknowledgments}}. It is, of course, impossible to thank here all the people whose work had a strong influence on the research presented here. I would nevertheless like to express my gratitude to the following people for their help throughout the years ~: P. Baumann, F. Bergeron, R. Bezrukavnikov, T. Bozec, I. Burban, P.-H. Chaudouard, I. Cherednik, W. Crawley-Boevey, B. Davison, B. Enriquez, P. Etingof, B. Feigin, D. Fratila, I.~Frenkel, T. Hausel, J. Heinloth, D. Hernandez, S.-J. Kang, M. Kapranov, M.~Kashiwara, B. Keller, M. Khovanov, G. Laumon, B. Leclerc, E. Letellier, G. Lusztig, D. Maulik, S. Meinhardt, A. Mellit, A.~Minets, S. Mozgovoy, H. Nakajima, A. Negut, A. Okounkov, S. Peng, P.-G. Plamondon, S. Riche, C. Ringel, F. Rodriguez-Villegas, R. Rouquier, M. Rosso, F. Sala, P.~Samuelson, Y. Soibelman, M. Reineke, M. Varagnolo, E. Vasserot,  B. Webster and J. Xiao.

\end{document}